%% file: main.tex
\documentclass[10pt]{report}
\usepackage
[
        a4paper,
        left=0.9cm,
        right=0.9cm,
        top=2cm,
        bottom=2cm,
]
{geometry}
\usepackage[T1]{fontenc}
\usepackage[latin9]{inputenc}
\usepackage{pdfpages}

\input{includes}
\input{defs}

\newboolean{ThesisIntro}
\setboolean{ThesisIntro}{true}

\begin{document}
\onehalfspacing

\begin{titlepage}
	\centering
	{\scshape\Large 
	The Hebrew University of Jerusalem \\ the Faculty of Mathematics and Natural Science \\ Einstein Institute of Mathematics
	\par}
	\vspace{1cm}
	{\Large
	In partial fulfillment of the requirements for the degree of
	\par}
	\vspace{0.5cm}
	{\Large
	Master of Science in Mathematics
	\par}
	\vspace{1.5cm}
	{\huge\bfseries
	Geometric branching operators on the Bruhat-Tits building associated to the symplectic group
	\par}
	\vspace{2cm}
	{\Large\itshape
	Submitted by Tal Nelken
	\par}
	\vfill
	Supervised by \par
	Prof. Ori Parzanchevski
	\vfill
	{\large \today \par}
\end{titlepage}


\input{Abstract}

\pagebreak

\tableofcontents
\chapter{Introduction}\label{sec:introintroduction}
\section{Background}
\input{Chapter1/Background_and_Introduction}

\section{Main results}
\input{Chapter1/mainResults}

\section{Outline of the thesis}\label{subsec:paper outline}
\input{Chapter1/OutlineOfPaper}

\chapter{Preliminaries and Notations}\label{sec:Preliminaries}

\section{Local fields and p-adic number fields}
\input{Chapter2/IntroFields}

\section{The Symplectic Group}
\input{Chapter2/IntroSymplectic}

\section{Geometric branching operators}
\input{Chapter2/IntroOperators}
\label{IntroOperators}

\section{Lattices}
\input{Chapter2/IntroLattices}
\label{IntroLattices}

\section{Weyl Groups}
\input{Chapter2/IntroWeyl}

\section{Bruhat-Tits buildings}
\input{Chapter2/IntroBuildings}

\chapter{ The Bruhat-Tits building associated to $Sp_{2n}$}\label{sec:GSpn}
\section{Building the building}
\input{Chapter3/BuildingTheBuilding}\label{sec:BTB}
\subsection{Describing the building using the \# operator}\label{hashtag operator}
\input{Chapter3/BTBOperator}
\subsection{Types of vertices in the building}
\input{Chapter3/TypesVertices}
\section{The Weyl group and apartments}
\input{Chapter3/ApartmentWeyl}
\subsection{Cartan decomposition}
\input{Chapter3/IntroCartan}

\chapter{Geometric branching operators on $\Delta_{2}$}\label{sec:Main proof n=2}
\input{Chapter4/IntroGeoOp2}
\section{Describing the operators}
\input{Chapter4/DescribeOperators2}
\section{Proof that the SRW can be described by the operators}
\label{SRW2}
\input{Chapter4/SRW2}

\chapter{Geometric branching operators on $\Delta_{n}$}\label{sec:Main proof gen n}
\input{Chapter5/IntroGeoOpn}
\section{Describing the operators}
\input{Chapter5/DescribeOperatorsn}
\label{descOpn}
\section{Proof that the SRW can be described by the operators}
\label{SRWn}
\input{Chapter5/SRWn}

\chapter{Discussion and open problems}
\input{Discussion}

\bibliography{Bib}
\bibliographystyle{plain}
\end{document}

%% file: includes.tex
\usepackage{titling}
\usepackage{color}
\usepackage{babel}
\usepackage{refstyle}
\usepackage{float}
\usepackage{amsmath}
\usepackage{amsthm}
\usepackage{makecell}
\usepackage{bbm}
\usepackage{amssymb}
\usepackage{stmaryrd}
\usepackage{setspace}
\usepackage{comment}
\usepackage{mathtools}
\usepackage{subcaption}
\usepackage{caption}
\usepackage{tikz}
\usepackage{rotating} 
\usepackage{fontawesome} 
\usepackage[misc]{ifsym} 
\usepackage{ifthen}
\usepackage{graphicx}
\usepackage{tikz-3dplot}
\usepackage{xstring}
\usepackage{xcolor}
\usepackage{soul}
\usepackage[colorinlistoftodos]{todonotes}
\usepackage[unicode=true,pdfusetitle,bookmarks=true,bookmarksnumbered=false,bookmarksopen=false,breaklinks=false,pdfborder={0 0 1},backref=false,colorlinks=true,linkcolor=blue]{hyperref}
\usepackage[nameinlink]{cleveref}
\usepackage{fdsymbol}

\usetikzlibrary{calc}
\usetikzlibrary{matrix}
\pgfmathsetmacro{\scaling}{0.4}

\usepackage{mleftright}
\mleftright 
\usepackage{etoolbox}
\usepackage{dynkin-diagrams}

%% file: defs.tex
\usetikzlibrary{patterns}
\tikzset{
        hatch distance/.store in=\hatchdistance,
        hatch distance=10pt,
        hatch thickness/.store in=\hatchthickness,
        hatch thickness=1pt
}
\makeatletter
\pgfdeclarepatternformonly[\hatchthickness]{my horizontal lines}
{\pgfpointorigin}{\pgfqpoint{100pt}{1pt}}{\pgfqpoint{100pt}{3pt}}
{
    \pgfsetcolor{\tikz@pattern@color}
    \pgfsetlinewidth{\hatchthickness}
    \pgfpathmoveto{\pgfqpoint{0pt}{0.5pt}}
    \pgfpathlineto{\pgfqpoint{100pt}{0.5pt}}
    \pgfusepath{stroke}
}

\pgfdeclarepatternformonly[\hatchthickness]{my vertical lines}
{\pgfpointorigin}{\pgfqpoint{1pt}{100pt}}{\pgfqpoint{3pt}{100pt}}
{
    \pgfsetcolor{\tikz@pattern@color}
    \pgfsetlinewidth{\hatchthickness}
    \pgfpathmoveto{\pgfqpoint{0.5pt}{0pt}}
    \pgfpathlineto{\pgfqpoint{0.5pt}{100pt}}
    \pgfusepath{stroke}
}

\pgfdeclarepatternformonly[\hatchthickness]{my north east lines}
{\pgfqpoint{-1pt}{-1pt}}{\pgfqpoint{4pt}{4pt}}{\pgfqpoint{3pt}{3pt}}
{
    \pgfsetcolor{\tikz@pattern@color}
    \pgfsetlinewidth{1pt}
    \pgfpathmoveto{\pgfqpoint{0pt}{0pt}}
    \pgfpathlineto{\pgfqpoint{3.1pt}{3.1pt}}
    \pgfusepath{stroke}
}

\pgfdeclarepatternformonly[\hatchthickness]{my north west lines}
{\pgfqpoint{-1pt}{-1pt}}{\pgfqpoint{4pt}{4pt}}{\pgfqpoint{3pt}{3pt}}
{
    \pgfsetcolor{\tikz@pattern@color}
    \pgfsetlinewidth{1.2pt}
    \pgfpathmoveto{\pgfqpoint{0pt}{3pt}}
    \pgfpathlineto{\pgfqpoint{3.1pt}{-0.1pt}}
    \pgfusepath{stroke}
}

\colorlet{grayX}{black!26!white}
\colorlet{darkgrayX}{black!47!white}
\makeatother

\makeatletter

\AtBeginDocument{}
\AtBeginDocument{}
\AtBeginDocument{}
\AtBeginDocument{}

\theoremstyle{definition}
\newtheorem{definition}{Definition}[section]
\newtheorem{remark}[definition]{Remark}
\newtheorem{claim}[definition]{Proposition}
\newtheorem{theorem}{Theorem}[section]
\newtheorem{corollary}[theorem]{Corollary}
\theoremstyle{plain}

\newtheorem*{example*}{Example}
\theoremstyle{plain}
\newtheorem*{thm*}{\protect\theoremname}
\theoremstyle{plain}

\theoremstyle{plain}

\theoremstyle{plain}

\theoremstyle{plain}

\theoremstyle{plain}

\crefname{equation}{}{}

\makeatother

\providecommand{\corollaryname}{Corollary}

\providecommand{\lemmaname}{Lemma}
\providecommand{\probname}{Problem}

\providecommand{\propositionname}{Proposition}

\providecommand{\theoremname}{Theorem}
\providecommand{\observationname}{Observation}

\makeatletter

\newif\if@mainmatter \@mainmattertrue
\newcommand\frontmatter{%
    \cleardoublepage
  \@mainmatterfalse
  \pagenumbering{Roman}}
\newcommand\mainmatter{%
    \cleardoublepage
  \@mainmattertrue
  \pagenumbering{arabic}}
\makeatother

\def\wasyfamily{\fontencoding{U}\fontfamily{wasy}\selectfont}
\DeclareTextFontCommand{\textwasy}{\wasyfamily}
\DeclareSymbolFont{wasy}{U}{wasy}{m}{n}
\SetSymbolFont{wasy}{bold}{U}{wasy}{b}{n}

\newcommand{\F}{\mathbb{F}}
\newcommand{\N}{\mathbb{N}}
\newcommand{\Z}{\mathbb{Z}}
\newcommand{\Q}{\mathbb{Q}}
\newcommand{\R}{\mathbb{R}}
\newcommand{\Qp}{\mathbb{Q}_{p}}
\newcommand{\cA}{\mathcal{A}}

\newcommand{\cW}{\mathcal{W}}
\newcommand{\Op}{\mathcal{O}}



%% file: Abstract.tex
\begin{abstract}
In this thesis we study geometric branching operators on simplices in the (affine) Bruhat-Tits building associated to the group $Sp_{n}(\F)$, for a local field $\F$. We use these operators to study the simple random walk on the building's vertices. This work lays the foundations for studying the total-variation cutoff phenomenon for simple random walk on Ramanujan complexes of type $\Tilde{C}_n$ ($n\geq 2$), as done in \cite{CP22} for Ramanujan complexes of type $\Tilde{A}_n$ ($n\geq 1$).
\end{abstract}

%% file: Chapter1/Background_and_Introduction.tex
A $d$-regular graph $X$ with $n$ vertices ($n,d \in \N$) is called an $\epsilon$-expander graph ($0 \less \epsilon \in \R$) if for every subset $A \subseteq X$, $|\partial A| \geq \epsilon (1-\frac{|A|}{|X|})|A|$, where $\partial A$, the boundary of $A$, is the set of all the vertices in $X$ of distance $1$ from the set $A$. 
A family of expanders is a family of $d$-regular graphs which are all $\epsilon$-expander graphs, for a global $\epsilon$.
Expander graphs are useful in mathematics and computer science because they are sparse graphs with strong connectivity properties. One of their properties is that their non-trivial adjacency spectrum is tightly bound. 
Expander graphs have been used, for example, in solving problems in the fields of error correcting codes, randomized polynomial-time algorithms, communication networks and more \cite{expanders}. 

A $d$-regular connected graph is called a Ramanujan graph if for every eigenvalue $\lambda$ of its adjacency matrix, $|\lambda|=d$ or  $|\lambda|\leq 2\sqrt{d-1}$. By the Alon-Boppana theorem \cite{expanders}, Ramanujan graphs are optimal expanders, from the spectral point of view. 

Ramanujan complexes are generalizations of Ramanujan graphs in higher dimensions.
They are defined as following: let $G$ be a simple algebraic group of rank $d$ over a non-archimedean local field, and let $B$ be the Bruhat-Tits (affine) building associated wih $G$. Let $C$ be a chamber in $B$ and $I$ the stabilizer of $C$ under the $G$-action on $B$. Let $\Gamma$ be a torsion-free lattice in $G$. 
Then $X=\Gamma \backslash B$ is a finite $d$-dimensional simplicial complex. $X$ is called a Ramanujan complex if every irreducible infinite-dimensional $G$-subrepresentation of $L^2(\Gamma \backslash G)$ containing an $I$-fixed vector, is weakly contained in $L^2(G)$.

Simple random walks (SRWs) on Ramanujan graphs exhibit the cutoff phenomenon, namely for a relatively long period of time the distribution of the SRW is almost as far from the stationary distribution as possible, then over a short period of time (the cutoff window) it mixes completely and becomes very close to the stationary distribution \cite{LP16}. This assertion invites the question, do random walks on Ramanujan complexes exhibit cutoff?
In \cite{CP22}, cutoff was established for SRW on Ramanujan complexes arising from the group $PGL_n(\F)$, for a non-archimedean local field $\F$. The cutoff was established using the following method: the SRW was factorized using certain "geodesic" geometric branching operators on the building's cells, which were studied in \cite{LLP20}.

This research has two purposes: the first purpose is to study geometric branching operators on cells of the building associated to $Sp_n(\F)$, for a non-archimedean local field $\F$. The second purpose is to use these operators in order to study the SRW on Ramanujan complexes arising from this group.

%% file: Chapter1/mainResults.tex
First, geometric branching collision-free regular operators (see \ref{IntroOperators} for definitions) were found on the chambers of the Bruhat-Tits (affine) building $\Delta_n$ associated to the group $GSp_{2n}(\F)$ for a non-archimedean local field $\F$. Next, by defining operators from vertices to chambers in the building and vice versa, these operators were combined to create geometric branching operators on the vertices of $\Delta_n$. Last, a section of the building was identified with a subset of $\N^n$ and it was shown all the vertices in the section are in the image of the above operators.

This research suggests a way to generalize the process of finding geometric branching operators on the vertices of Bruhat-Tits affine buildings associated to a classical group over a local field.
The process would be as following: first, find the translation elements corresponding to the coroot lattice in the building's associated affine Weyl group $\mathcal{W}$. 
Given a translation element $t\in \mathcal{W}$, we can define a branching operator $T$ on the building's chambers by $T(C)=ItC$, where $C$ is a chamber and $I$ is its stabilizer in $G$. This operator is geometric, and if the group $G$ is a simple group, it is also collision-free. It will also be $k$-regular for $k=|ItC|$.
Thus, this process will enable the construction of a family of geometric branching collision-free and regular operators on the building's chambers.
Second, find the orbits of the $G$-action on the building's vertices. For every orbit, define an operator from the set of chambers to that orbit using stabilizers, like the operators $T_{n,0}^{(i)}$ in \ref{descOpn}. 
Third, define an operator from a chosen special vertex $\xi$ to the chambers in the building, like $T_{0,n}$ in \ref{descOpn}. 
Last, combine these operators to give an operator from the orbit of $\xi$ to the set of all vertices in the building, like in \ref{SRW2} and \ref{SRWn}.

%% file: Chapter1/OutlineOfPaper.tex
\Cref{sec:Preliminaries} gives a brief taste of the main definitions and theoretical background required for this research.  \Cref{sec:GSpn} describes the Bruhat-Tits (affine) building associated to $GSp_{2n}(\F)$ for a local field $\F$. \Cref{sec:Main proof n=2} is dedicated to the results for the two-dimensional case, while \Cref{sec:Main proof gen n} describes the results in the general case.  

%% file: Chapter2/IntroFields.tex
Throughout this thesis, we will be working over local fields (see definition \ref{loclf}). A useful example for such a field is the p-adic field $\Q_{p}$ for a prime number $p\in\Z$, which I will briefly describe in this section, following \cite{padics}. 

Fix a prime number $p\in\Z$.
\begin{definition}[$p$-adic valuation]
The $p$-adic valuation on $\Z$ is the function $\nu_{p}:\Z-\{0\}\longrightarrow\Z$ which is defined as follows: let $0\neq n\in\Z$. $\nu_{p}(n)$ is the unique non-negative integer satisfying 
$n=p^{\nu_{p}(n)}n' \quad with \quad p\nmid n'.$

We extend $\nu_{p}$ to $\Q$ in the following way: if $x=a/b\in\Q^{\times}$, then $$\nu_{p}(X)=\nu_{p}(a)-\nu_{p}(b),$$ and we set $\nu_{p}(0)=\infty$.
\end{definition}

In fact, the $p$-adic valuation on any $x\in\Q^{\times}$ satisfies the following equation: $$x=p^{\nu_{p}(X)}\cdot\frac{a}{b} \quad p\nmid ab.$$

\begin{definition}[$p$-adic absolute value]
The $p$-adic absolute value $|\cdot|_{p}$ on $\Q$ is defined as follows: for $0\neq x \in \Q,$ $$|x|_{p}=p^{-\nu_{p}(x)}.$$
For $x=0$ we define $|0|_{p}=0.$
\end{definition}

It is easy to check that the $p$-adic absolute value is a non-archimedean absolute value on $\Q$. Thus, it is possible to define the $p$-adic metric on $\Q$ in the following way: $$d_{p}:\Q\times\Q\longrightarrow\R_{\geq0} \quad \quad d_{p}(a,b)=|a-b|_{p}.$$

\begin{definition}[$\Q_{p}$]
The $p$-adic field $\Q_{p}$ is defined to be the completion of $\Q$ with respect to the $p$-adic metric.
\end{definition}

The above definition is useful but it's difficult to deduce from it how $p$-adic numbers actually look. The next claim gives a concrete description of $p$-adic numbers.

\begin{claim}[Elements in $\Q_{p}$ are power series in p]
Every $x\in\Q_{p}$ can be written as $$x=a_{-m}p^{-m}+\cdot\cdot\cdot+a_{-1}p^{-1}+a_{0}+a_{1}p+a_{2}p^{2}+\cdot\cdot\cdot+a_{n}p^{n}+\cdot\cdot\cdot=\sum_{n\geq-m} a_{n}p^{n}$$ with $0\leq a_{n}\leq p-1$ and $-m=\nu_{p}(x)$ if $a_{-m} \neq 0$. Additionally, this representation is unique.
\end{claim}

\begin{claim} [The $p$-adic valuation extends to $\Qp$]
For every $0\neq x\in\Qp$ there exists an integer $\nu_{p}(x)$ such that $|x|_{p}=p^{-\nu_{p}(x)}$, and we define $\nu_{p}(0)=\infty$. The resulting function $\nu_{p}$ is a discrete valuation on $\Qp$.
\end{claim}

\begin{definition}[$p$-adic integers]
The ring of $p$-adic integers is defined to be $$\Z_{p}=\{x\in\Qp:|x|_{p}\leq 1\}.$$ 
\end{definition}

\begin{claim}[$p$-adic residue field]
The ring $\Z_{p}$ has a maximal ideal which is $p\Z_{p}=\{x\in\Qp:|x|_{p}\less 1\}$.Thus the $p$-adic residue field is the finite field with p elements $\F_{p}$.
\end{claim}

\begin{definition}[local field]
\label{loclf}
A field $\F$ is called a local field if it is a discrete valuation field (with non-trivial valuation) which is complete with respect to the topology induced by the valuation, and it has a finite residue field.
\end{definition}

\begin{remark}
The $p$-adic valuation on the field $\Qp$ is a discrete valuation thus $\Qp$ is a discrete valuation field with a finite residue field. Furthermore, $\Qp$ is complete with respect to the topology induced by the $p$-adic valuation $\nu_{p}$. Thus, $\Qp$ is a local field.
\end{remark}

%% file: Chapter2/IntroSymplectic.tex
The goal of the following section is to give a brief introduction to the symplectic group and the group of symplectic similitudes over a general field in the context needed for this thesis, following \cite{Symplectic}. Throughout this section, let V be a vector space of finite dimension over a field $\F$, and let $B:V\times V\longrightarrow \F$ be a bilinear form.

\begin{definition}[Anti-symmetric bilinear form]
A bilinear form $B$ is called anti-symmetric (or skew-symmetric) if for every $v,w\in V$ $B(v,w)=-B(w,v)$.
\end{definition}

Suppose from now on that $B$ is non-degenerate and anti-symmetric. The following proposition summarizes some theorems in chapters 2 and 3 in \cite{Symplectic}.

\begin{claim}
The matrix form of $B$ is congruent to $\Omega=\big(\begin{smallmatrix}
0 & I \\
-I & 0
\end{smallmatrix}\big)$. In particular, V has even dimension over $\F$ and a basis of mutually orthogonal hyperbolic pairs. 
\end{claim}

\begin{definition}[The symplectic group]
The symplectic group, $Sp_{2n}(\F)$, is defined to be the group of automorphisms of $\Omega$, $Aut(\Omega)=\{A\in GL_{2n}(\F) : A^{t}\Omega A = \Omega\}$.
\end{definition}

\begin{claim}
For every $A\in Sp_{2n}(\F), \quad det(A)=1$.
\end{claim}
\begin{proof}
$Sp_{2n}(\F)$ is generated by the symplectic transvections, which are elements in $SL_{2n}(\F)$ \cite{Symplectic}. Thus $Sp_{2n}(\F) \leq SL_{2n}(\F)$. 
\end{proof}

\begin{definition}[The group of symplectic similitudes]
The group of symplectic similitudes, $GSp_{2n}(\F)$, is defined to be the group $\{A\in GL_{2n}(\F) : A^{t}\Omega A = \mu(A)\Omega\}$ such that $\mu(A)\in\F^{\times}$.
\end{definition}

\begin{definition} 
\label{symplectic transversion}
We define for every $A\in GSp_{2n}(\F)$ $A^{\#} = -\mu(A)\Omega A^{t^{-1}}\Omega$. 
\end{definition}

\begin{remark}
$GSp_{2n}(\F)$ is the set of fixed points under the operator $\#$, for matrices in $GL_{2n}(\F)$ for which there exists a number $\mu(A) \in \F^{\times}$, 
because if $A\in GSp_{2n}(\F)$, then $$A^{t}\Omega A = \mu(A)\Omega$$ and A is invertible. 
This is equivalent to 
 $$-\mu(A)\Omega A^{t^{-1}}\Omega =A.$$ Thus, $A\in GSp_{2n}(\F) \Longleftrightarrow A^{\#} = A$.
\end{remark}

Now we assume $\F$ is a local field with a uniformizing parameter $\varpi$ and discrete valuation $\nu$. The following claim is an easy exercise and will prove to be useful in the following chapters:

\begin{claim}
\label{val and det in GSp}
Let $A\in GSp_{2n}(\F)$.
\begin{enumerate}
    \item $det(A)^{2} = \mu(A)^{2n}$
    \item $n\mid \nu(det(A))$.
    \item $\nu(\mu(A))\equiv \frac{\nu(det(A))}{n} \mod{2}$
\end{enumerate}
\end{claim}

%% file: Chapter2/IntroOperators.tex
In this section we define operators which we can think about as describing the allowed moves in a random walk.

\begin{definition} [branching operator]
Let $X$ be a set. A branching operator on $X$ is a map $T:X\longrightarrow2^{X}$.
\end{definition}

\begin{remark}
\begin{enumerate}
    \item For a branching operator $T$ describing the allowed moves of a random walk, we can describe the random walk by starting with a point $x\in X$, and in the next step go to some uniformly chosen $y\in T(x) \subseteq X$.
    \item $T$ describes and can be described by the digraph $(X,\{(x\longrightarrow y)|x\in X, y\in T(x)\})$.
\end{enumerate}
\end{remark}

\begin{definition} [geometric branching operator]
Let $G$ be a group and let $X$ be a $G$-set. Let $T$ be a branching operator on $X$. $T$ is geometric if it satisfies $\forall x \in X, \forall g\in G, T(gx)=gT(x).$
\end{definition}

Let $T$ be a geometric branching operator on $X$. Assume that the $G$-action on $X$ is transitive. Let $x_{0}\in X$ and let $K=Stab_{G}(x_{0})$, thus by the orbit-stabilizer theorem, $X\cong G/K$ as $G$-sets. 

\begin{claim}
\label{geoOp}
\begin{enumerate}
    \item $T(x_{0})$ determines $T$.
    \item $T(x_{0})$ is $K$-stable.
    \item Any $K$-stable set $Y\subseteq X$ defines a geometric branching operator, via $T(x_0)=Y$.
\end{enumerate}
\end{claim}
\begin{proof}
\begin{enumerate}
    \item Let $x\in X$, then there exists a $g\in G$ such that $x=gx_{0}$. Thus, $T(x)=T(gx_{0})=gT(x_{0})$.
    \item Let $k\in K$. $kT(x_{0})=T(kx_{0})=T(x_{0})$. Thus $KT(x_{0})=T(x_{0})$.
    \item Define $T(x_{0})=Y$. From 1., this defines $T$. And $T$ is well-defined, because if $gx_{0}=g'x_{0}$, then $g^{-1}g'\in K$, Thus $g^{-1}g'T(x_{0})=g^{-1}g'Y=Y=T(x_{0})\Longrightarrow T(gx_{0})=gT(x_{0})=g'T(x_{0})=T(g'x_{0})$.
\end{enumerate}
\end{proof}

Let $s\in G$ and $s_{1},...s_{r}\in G$ such that $KsK=\bigsqcup_{i=1}^{r}s_{i}K$. The set $\{s_{i}x_{0}\}_{i=1}^{r}$ is $K$-stable, as it is equal to the set $Ksx_{0}$ (because $Ksx_{0}=KsKx_{0}=\bigsqcup_{i=1}^{r}s_{i}Kx_{0}=\{s_{i}x_{0}\}_{i=1}^{r}$). Thus, by the previous claim, we can define a geometric branching operator $T_{s}$ by $T_{s}(x_{0})=\{s_{i}x_{0}\}_{i=1}^{r}$.

\begin{claim}
Assume that every double coset in $K\backslash G / K$ is a finite union of cosets in $G/K$.
\begin{enumerate}
    \item If $\emptyset \neq Y\subseteq X$ is finite and $K$-stable, then $Y$ is of the form $KSx_{0}$ for some finite $S\subseteq G$.
    \item Any geometric branching operator on $X$ such that $T(x_{0})$ is finite is of the form $T(x)=T_{s_{1}}(x)\bigcup...\bigcup T_{s_{r}}(x)$ for $s_{1},...s_{r}\in G$.
\end{enumerate}
\end{claim}

\begin{remark}
In this work it will often be the case that double cosets are a finite union of left cosets, as assumed in this proposition.
\end{remark}
\begin{proof}
\begin{enumerate}
    \item $Y \neq \emptyset$ so there exists some $x \in Y$, and denote it by $x_0$. $G$ acts on $X$ transitively so for every $y \in Y$ there exists (at least one) $g_y$ such that $g_yx_0=y$. Choose one such $g_y$ for every $y\in Y$ and define the set $S$ to be $S=\{g_y|y\in Y\}$. $S$ is finite because $Y$ is finite. Thus, $Sx_0=Y \Rightarrow KSx_0=KY=Y$ because $Y$ is $K$-stable.
    \item $T(x_{0})=Y$ is finite, so from the first part of this claim, there exists a finite $S\subseteq G$ such that $T(x_{0})=KSx_{0}.$ Denote $S=\{s_{1},...,s_{r}\}$. Thus $T(x_{0})=Ks_{1}x_{0} \bigcup ... \bigcup Ks_{r}x_{0} = T_{s_{1}}(x_{0})\bigcup...\bigcup T_{s_{r}}(x_{0}) \Longrightarrow T(x)=T_{s_{1}}(x)\bigcup...\bigcup T_{s_{r}}(x)$.
\end{enumerate}
\end{proof}

Now let's assume that the $G$-action on $X$ is not transitive, but has $n$ orbits $Gx_{1},...,Gx_{n}$. Thus $X=\bigsqcup_{i=1}^{n} Gx_{i}$. In this case $T$ is determined by $T(x_{1}),... ,T(x_{n})$ as following: 
Denote the stabilizer of $x_{i}$ by $K_{i}$. There exist subsets $S_{i,j}\subseteq G$ for $1\leq i,j \leq n$ such that $$T(x_{i})=S_{i,1}x_{1}\bigsqcup S_{i,2}x_{2}\bigsqcup...\bigsqcup S_{i,n}x_{n},$$ and $T(x_{i})$ is $K_{i}$-stable. 
In particular, $K_iS_{i,j}K_j=S_{i,j}K_j$ for all $1\leq i,j \leq n$.

We give two last definitions for branching operators:

\begin{definition} [collision-free operator]
We call a branching operator on $X$ collision-free if its associated digraph has at most one path from $x$ to $y$ for every $x,y\in X$.
\end{definition}

\begin{definition} [regular branching operator]
A $k$-branching operator $T$ on $X$ is a map from $X$ to all the $k$-element subsets of $X$, or equivalently the associated digraph is $k$-out regular. $T$ is called a $k$-regular branching operator if the associated digraph is both $k$-out and $k$-in regular.
\end{definition}

\begin{remark}
In the paper \cite{LLP20}, the main result states that if $T$ is a $k$-regular geometric collision-free branching operator on a Bruhat-Tits building $\mathcal{B}$, the associated random walk on an $n$-vertex Ramanujan complex with universal cover $\mathcal{B}$ has cutoff at time $log_{k}n$.
\end{remark}

%% file: Chapter2/IntroLattices.tex

In this section we will give a brief introduction to lattices in vector spaces. We will use lattices as representatives for vertices in buildings in the following chapters. Throughout this chapter let $\F$ be a local field and $\mathcal{O}$ its ring of integers. Let $V$ be a finite-dimensional vector space over $\F$. Let $\varpi$ be a uniformizing parameter, and $F = \mathcal{O} / \varpi \mathcal{O}$ be the residue field of $\F$.

\begin{definition} [$\mathcal{O}$-lattice]
A $\mathcal{O}$-lattice $L$ of $V$ is a free $\mathcal{O}$-submodule of $V$ of rank $dim_{\F}(V)$.
\end{definition}

\begin{example*}
Let $n=dim_{\F}(V)$, and $v_1,...,v_n$ a $\F$-basis for $V$. Then $L=\mathcal{O}v_1 + \mathcal{O}v_2 + ... + \mathcal{O}v_n$ is a $\mathcal{O}$-lattice. In fact, any $\mathcal{O}$-lattice is in this form for some basis of $V$.
\end{example*}

\begin{definition}[homothetic lattices]
Let $L$ and $L'$ be $\mathcal{O}$-lattices of $V$. $L$ and $L'$ are called homothetic if there exists $\alpha \in \F^{\times}$ such that $L'=\alpha L$. 
\end{definition}

\begin{remark}
Homothety of $\mathcal{O}$-lattices is clearly an equivalence relation. Denote by $[L]$ the homothety class of a $\mathcal{O}$-lattice $L$, i.e. the equivalence class of $L$ under the homothety relation.
\end{remark}

\begin{definition}[primitive $\mathcal{O}$-lattice]
Let $\langle \cdot , \cdot \rangle$ be a non-degenerate alternating bilinear form on $V$ (as in a symplectic vector space, for example). Call a $\mathcal{O}$-lattice L primitive if $\langle L , L \rangle \subseteq \mathcal{O}$, and $\langle \cdot , \cdot \rangle$ induces a non-degenerate alternating form on the vector space $L/\varpi L$ over $F$.
\end{definition}
\label{primLatt}

%% file: Chapter2/IntroWeyl.tex
The current and following chapters are on subjects which can fill entire books. Accordingly, these chapters will only "touch the tip of the iceberg" on these subjects. In both chapters we first follow \cite{buildings3}, then we follow \cite{buildings2}. let $V$ be a Euclidean space, i.e. a real finite-dimensional vector space with an inner product $\langle \cdot , \cdot \rangle $.

\begin{definition}[hyperplane]
A hyperplane $H\subseteq V$ is a codimension $1$ subspace of $V$.
\end{definition}

\begin{definition}[reflection along a hyperplane]
Let $H$ be a hyperplane. A reflection along $H$ is a linear transformation $s_{H}:V\longrightarrow V$ which is the identity on $H$ and is the multiplication by $-1$ on the orthogonal complement of $H$.
\end{definition}

\begin{definition}[finite reflection group]
Let $\mathcal{H}$ be a set of hyperplanes. A finite $\mathcal{H}$-reflection group is a finite group $W\leq GL(V)$ which is generated by reflections $s_{H}$ such that $H\in \mathcal{H}$.
\end{definition}

\begin{claim}
\label{finiteRefG}
Let $R$ be a finite set of non-zero vectors in $V$. Let $\mathcal{H}=\{\alpha^{\bot} | \alpha \in R\}$ be a set of hyperplanes and $W=\langle s_{H} | H \in \mathcal{H} \rangle$. If $R$ is invariant under the action of $W$, then $W$ is finite.
\end{claim}

\begin{remark}
Let $\mathcal{H}$ be a set of hyperplanes such that $W$ defined as above is a finite reflection group. The hyperplanes in $\mathcal{H}$ cut $V$ into polyhedral sections which are cones over simplices, called cells. 
Thus we obtain a simplicial complex, whose intersection with the unit sphere in $V$ triangulates it. 
A chamber is an open convex set which is a connected component of the complement $V-\bigcup_{H\in \mathcal{H}}H$. Define $\Sigma$ to be the poset of closed cells, ordered by inclusion.
\label{chamber}
\end{remark}

From now on let $\mathcal{H}$ be a set of hyperplanes such that $W$ generated by reflections along the hyperplanes in $\mathcal{H}$ is a finite reflection group. Assume that $\mathcal{H}$ is $W$-invariant. Note that $W$ acts on $\Sigma$. 

\begin{theorem}
\begin{enumerate}
    \item Let $C$ be a chamber. $W$ is generated by the reflections $s_{H}$ such that $H$ is a wall (a face of codimension 1) of $C$.
    \item Let $H_1,...,H_r$ be the set of walls of a chamber $C$, with normals $e_1,...,e_r$ such that the normals point into $C$. Let $s_1,...,s_r$ be the reflections along the walls. Assume that $dimV=r$. Then $\langle e_i,e_j\rangle =-cos(\frac{\pi}{m_{ij}})$ for $1\leq i,j \leq r $ and $m_{ij}$ the order of $s_is_j$.
\end{enumerate}
\end{theorem}

\begin{corollary}
In the conditions of the above theorem, $W$ is completely determined, up to isomorphism, by the $r\times r$ matrix $M=(m_{ij})$. $M$ is called the Coxeter matrix associated to $W$.
\end{corollary}

\begin{definition}[crystallographic root system]
A crystallographic root system $\Phi \subseteq V$ is a finite set of non-zero vectors which satisfy the following:
\begin{enumerate}
    \item $\forall \alpha \in \Phi, \quad \Phi \cap \R \alpha=\{ \alpha, -\alpha \}$
    \item $\forall \alpha \in \Phi, \quad s_{\alpha}(\Phi)=\Phi$
    \item $\forall \alpha,\beta \in \Phi, \quad \frac{2\langle \alpha , \beta \rangle}{\langle \alpha , \alpha \rangle} \in \Z$
\end{enumerate}
\end{definition}

\begin{definition}[irreducible crystallographic root system]
A crystallographic root system $\Phi$ is called irreducible if $\forall \alpha, \beta \in \Phi,$ $\exists \alpha=\alpha_1,...,\alpha_m,\alpha_{m+1}=\beta \in \Phi$ such that $\prod_{i=1}^{m}\langle \alpha_{i},\alpha_{i+1} \rangle \neq 0$.
\end{definition}

\begin{remark}
Many crystallographic root systems arise from representations of semi-simple Lie groups.
\end{remark}

\begin{definition}[Weyl group]
The (spherical) Weyl group of a crystallographic root system $\Phi$, is $W(\Phi)=\langle s_{\alpha} | \alpha \in \Phi \rangle.$
\end{definition}

\begin{remark}
By proposition \ref{finiteRefG}, the Weyl group of $\Phi$, for a crystallographic root system $\Phi$, is a finite reflection group, satisfying all the conditions of the above theorem.
\end{remark}

\begin{definition}
The pair $W,S$ such that $W$ is a reflection group generated by $S$ is called a Coxeter System.
\end{definition}

\begin{remark}
Coxeter groups are abstract reflection groups which can be seen as abstract generalizations of Weyl groups.
\end{remark}

\begin{definition}[Coxeter diagram]
The Coxeter diagram is a diagram encoding all the information in a Coxeter matrix $M$. The diagram is a graph with $n$ vertices, such that $n=|S|$ for the set of generators $S$ of $W$. Label the vertices with the indices $1,...,n$. Vertices $i$ and $j$ are connected by an edge if and only if $m_{ij}=order(s_{i}s_{j})\geq 3$. If $m_{i,j}\geq4$ the edge is labeled with $m_{ij}$. Coxeter diagrams can be seen as simplified versions of Dynkin diagrams, which have some additional information.
\end{definition}

\begin{remark}
There is a classification of all the irreducible root systems, up to isomorphism; Three infinite families: $A_n (n\geq1),B_n=C_n(n \geq 2),D_n(n\geq4)$, and some exceptional groups. The infinite families have the following Coxeter diagrams:

$A_n$ \dynkin[labels={1,2,n-1,n}, scale=2, Coxeter]A{}

$B_n=C_n$ \dynkin[labels={1,2,3,n-2,n-1,n}, scale=2.3,Coxeter]B{}

$D_n$ \dynkin[labels={1,2,n-2,n-1,n,n+1},label directions={,,,right,,}, scale=1.8,
Coxeter]D{}
\end{remark}

\begin{definition}[affine hyperplane]
Let $\alpha \in \Phi$ and $k\in \R$. An affine hyperplane $H_{\alpha,k}$ is the set $\{v \in V | \langle v,\alpha \rangle = k\}$. 
\end{definition}

\begin{definition}[coroot]
Let $\alpha \in \Phi$. $\alpha^{\vee}=\frac{2}{\langle \alpha , \alpha \rangle}\alpha$ is called the coroot of $\alpha$. Define $\Phi^{\vee}=\{\alpha^{\vee} | \alpha \in \Phi \}$ to be the dual root system.
\end{definition}

\begin{definition}[reflection along an affine hyperplane]
Let $\alpha \in \Phi$ and $k \in \R$, define $s_{\alpha,k}:V\longrightarrow V$ by $s_{\alpha,k}(v)=v-(\langle v,\alpha \rangle -k)\alpha^{\vee}.$ $s_{\alpha,k}$ is the reflection along the affine hyperplane $H_{\alpha,k}$.
\end{definition}

\begin{definition}[root and coroot latice]
$L(\Phi)=Span_{\Z}(\Phi)$ is called the root lattice. Accordingly, $L(\Phi^{\vee})$ is called the coroot lattice. These lattices act on $V$ by translations.
\end{definition}

\begin{definition}[affine Weyl group]
The affine Weyl group of $\Phi$ is 
$W_a(\Phi)=\langle s_{\alpha,k} | \alpha \in \Phi,k\in \Z \rangle$.
\end{definition}

\begin{example*}
\label{roots coroots C2}
Let $e_1,e_2$ be the standard basis vectors for $\R^2$ with the standard inner product. Define $\alpha_1=2e_1$, $\alpha_2=-e_1+e_2$. Then $\Phi=\{\alpha_1, \alpha_2, \alpha_1+\alpha_2, \alpha_1+2\alpha_2, -\alpha_1, -\alpha_2, -\alpha_1-\alpha_2, -\alpha_1-2\alpha_2\}=\{\pm 2e_1, \pm 2e_2, \pm e_1 \pm e_2\}$ is a crystallographic root system. 
It is of type $C_2$ and arises from the adjoint representation of the lie algebra of $Sp_4(\R)$. The coroots are $\Phi^{\vee}=\{\frac{\alpha_1}{2}, \alpha_2, \alpha_1+\alpha_2, \frac{\alpha_1+2\alpha_2}{2}, -\frac{\alpha_1}{2}, -\alpha_2, -\alpha_1-\alpha_2, \frac{-\alpha_1-2\alpha_2}{2}\}=\{\pm e_1, \pm e_2, \pm e_1 \pm e_2\}$.
\end{example*}

\begin{claim}
$W_a(\Phi)$ is the semi-direct product of $W(\Phi)$ with the (abelian) translation subgroup $L(\Phi^{\vee})$.
\end{claim}

Now, take $\mathcal{H}$ to be the set $H_{\alpha,k}$ for $\alpha \in \Phi$ and $k \in \Z$. Note that a chamber $C$ defined in \ref{chamber} is now a simplex. Similarly to the spherical case, $W_a(\Phi)$ is generated by the reflections across the faces of $C$. 
Additionally, the simplicial complex defined by the set of hyperplanes $\mathcal{H}$ is called the affine Coxeter complex associated to $W_a(\Phi)$. It can be seen as taking the chamber $C$ and repeatedly reflecting it by words in the generating reflections, giving a tessellation of $V$.

\begin{definition}[length function]
Let $S$ be a set of reflections across the faces of $C$, which generate $W_{a}(\Phi)$. Define the length function $l:W_{a}(\Phi)\longrightarrow \Z_{\geq0}$ by $$l(w)=min\{r|\exists s_1,...,s_r \in S; w=s_1\cdot \cdot \cdot s_r\}.$$ $l(w)$ is called the length of $w$, and any expression $w=s_1\cdot \cdot \cdot s_{l(w)}$ is called a reduced word for w.
\end{definition}

\begin{definition}[gallery, combinatorial distance]
Call chambers $C$ and $C'$ adjacent if they have a common codimension $1$ face.
A gallery from $C_0$ to $C_d$ is a sequence of chambers $\Gamma=(C_0,...,C_d)$ such that consecutive chambers are adjacent. $d$ is called the length of the gallery $\Gamma$. 
Define the combinatorial distance between chambers $C$ and $C'$ to be $$n(C,C')=min\{d|\exists \Gamma \quad of \quad length \quad d\quad s.t. \quad \Gamma \quad is \quad a \quad gallery \quad between \quad C \quad and \quad C' \}$$.
\end{definition}

For two chambers $C$ and $C'$, let $l(C,C')$ be the number of hyperplanes which separate between $C$ and $C'$.

\begin{claim}
Let $C$ and $C'$ be chambers, so that $wC=C'$. Then $l(C,C')=n(C,C')=l(w)$.
\end{claim}

%% file: Chapter2/IntroBuildings.tex
Affine Bruhat-Tits buildings have many definitions, but in general they can be thought of as very complicated simplicial complexes which are made up of many subcomplexes isomorphic to affine Coxeter complexes, called apartments, glued together in complex ways. 
Given an isometry group over a local field, an affine building can be associated to it by a construction using homothety classes of lattices over that field, and in turn the group acts on the simplices in the building, by acting on the homothety classes of the lattices. The building can also be constructed in an equivalent way using the $BN$-pair construction \cite{buildings2}.

\begin{definition}[building]
A building is a simplicial complex $\Delta$ which can be expressed as a union of subcomplexes $\Sigma$ called aparments such that:
\begin{enumerate}
    \item Each apartment $\Sigma$ is isomorphic to a Coxeter complex.
    \item For any two simplices $A,B\in \Delta$, there is an apartment containing both of them.
    \item if $\Sigma $ and $\Sigma'$ are two apartments containing $A$ and $B$, then there is an isomorphism $\Sigma \longrightarrow \Sigma'$ fixing $A$ and $B$ pointwise.
\end{enumerate}
A building is called thick if every codimension one face of a chamber $C$ is a face of at least three chambers.
An affine building is a thick building such that every apartment is an affine Coxeter complex.
\end{definition}

\begin{remark}
One-dimensional affine buildings are trees.
\end{remark}

Let $G$ be a group acting by simplicial-complex automorphisms on a building, such that it has a $BN$ pair. Let $C$ be a chamber and let $I$ be its stabilizer in $G$. Let $(W,S)$ be the Coxeter system arising from $G,B,N$. Then the following equations, which we will call "Bruhat relations" hold:
$$\forall s\in S, w\in W, IwI\cdot IsI=IwsI \quad if \quad  l(w) \less l(ws)$$
$$\forall s\in S, w\in W, IwI\cdot IsI=IwsI\sqcup IwI \quad if \quad  l(ws) \less l(w)$$
\label{BruhatRelations}

\begin{remark}
\label{Omega}
If $G$ has a subgroup $G_0$ which has a strict $BN$ pair, then the Weyl group of $G$ is a semi-direct product of a group $\Omega$ isomorphic to $G/G_0$ with the Weyl group $W_0$ of $G_0$. For every $\sigma \in \Omega$ and $w \in W_0$, 
$\sigma I w I = I \sigma w I$. Additionally, $\sigma I=I \sigma=I \sigma I$.
\end{remark}

\begin{corollary}
$\forall w_1,w_2\in W, \quad Iw_1I\cdot Iw_2I\supseteq Iw_1w_2I.$
\end{corollary}

%% file: Chapter3/BuildingTheBuilding.tex
To begin, let's establish some notations which will be used in the current and following chapters. Let $\F$ be a local field, $\mathcal{O}$ its ring of integers, $\varpi$ a uniformizing parameter, $F = \mathcal{O} / \varpi \mathcal{O}$ the residue field of $\F$. Let $V$ be a symplectic vector space over $\F$, with the anti-symmetric non-degenerate bilinear form $B$. 
Denote the dimension of $V$ over $\F$ by $2n$. Fix a basis of $V$ $\{u_{1},...,u_{n},w_{1},...,w_{n}\}$ so that $(u_{i},w_{i})$ form a hyperbolic plane under $B$, so in this basis $B$ has matrix form $\Omega$. 

The Bruhat-Tits building $\Delta_{n}$ for $Sp_{2n}(\F)$ is a simplicial complex of dimension $n$. It can be described in many equivalent ways, of which we will use two; the first way describes the building as a part of the building for $GL_{2n}(\F)$ which is stable under the symplectic $\#$ operation defined in \ref{symplectic transversion}, with additional simplices formed under that operation \cite{T79}.
The second way describes the vertices as homothety classes of lattices in $V$, and the simplices of higher dimension are defined by sub-lattice relations. The building itself is the resulting flag complex \cite{S05}. This thesis will mainly use the first way which will soon be described in detail. 

Following \cite{S05}, denote an $\mathcal{O}$-lattice of the shape $\mathcal{O}\varpi^{a_{1}}u_{1}\oplus\cdot\cdot\cdot\oplus\mathcal{O}\varpi^{a_{n}}u_{n}\oplus\mathcal{O}\varpi^{b_{1}}w_{1}\oplus\cdot\cdot\cdot\oplus\mathcal{O}\varpi^{b_{n}}w_{n}$ by the exponents $(a_{1},...,a_{n};b_{1},...b_{n})$, and its homothety class by $[a_{1},...,a_{n};b_{1},...b_{n}]$ \label{homothetyExponents}. 
Call a lattice symplectic if it has an $\mathcal{O}$-basis which is a symplectic similitude basis for $V$ with respect to $B$. Call a homothety class symplectic if it has a symplectic representative.
Associate a matrix $A$ to a general homothey class $[L]$ as following: take $L$ to be the unique representative of $[L]$ such that there exists a primitive lattice $\Gamma$, such that $\varpi \Gamma \less L \less \Gamma$, and $\langle L,L\rangle \subseteq \mathcal{O}$. Let $A$ be a matrix such that $A\mathcal{O}^{2n}=L$.
For a matrix $A\in M_{2n}(\mathcal{O})$, associate a lattice by taking the homothety class of the $\mathcal{O}$-span of $A$'s columns. Note that a homothety class $L$ is symplectic if and only if the associated matrix $A$ is in $GSp_{2n}(\mathcal{O})$.
Denote the apartment of $\Delta_{n}$ determined by $\{\F u_{i},\F w_{i}\}$ as $\cA$ and the Weyl group by $\cW$. $\cW$ is of type $\tilde {C}_{n}$, and it has a generator set $S$ where $|S|=n+1$ ($S$ generates the Weyl group associated to $Sp_{2n}(\F)$), and additional elements which aren't generated by $S$, which we will denote by $t$ (these are elements in $\Omega$ as in \ref{Omega}).

%% file: Chapter3/BTBOperator.tex
Denote the Bruhat-Tits building for $GL_{2n}(\F)$ by $\Theta_{2n}$. $\Theta_{2n}$ can be described as following: It's vertices are homothety classes of $\mathcal{O}$-lattices in $\F^{2n}$. Homothety classes $c_{0},...,c_{j}$ form a cell if, possibly after reordering, there are lattices $L_{i}\in c_{i}$ such that $\varpi L_{0}\less L_{j} \less L_{j-1} \less \cdot \cdot \cdot \less L_{0}$ \cite{buildings}. Thus, $\Theta_{2n}$ is an $(2n-1)$-dimensional simplicial complex. The set of vertices of $\Delta_{n}$ is the union of the set of vertices of $\Theta_{2n}$ which are stable under the $\#$ operation (where a lattice L is represented by the matrix associated to it), which are called 'special vertices', with a set of 'new' vertices - the center points of the edges between vertices of $\Theta_{2n}$ which are swapped under the $\#$ operation. These 'new' vertices formed by swapping two vertices from $\Theta_{2n}$ are the non-special vertices, and we will represent them by both original vertices joined by a $+$ sign. A cell is formed as in $\Theta_{2n}$ where for the non-special vertices we can, without loss of generality, choose one of the original homothety classes to be the representative of that vertex (as both those vertices were part of a simplex in $\Theta_{2n}$, we know that they have representatives forming a chain, so it doesn't matter which one we choose).  

\begin{example*}
\label{C2}
For $n=2$, the following homothety classes of lattices form a cell:
$$[0,0;0,0], \quad [0,1;0,0]+[1,1;1,0], \quad [1,1;0,0]$$
\end{example*}

In the following section we will see that $\Delta_{n}$ is an $n$-dimensional simplicial complex.

%% file: Chapter3/TypesVertices.tex
There are two types of vertices in $\Delta_{n}$ - special and non-special, as described in \ref{hashtag operator}. Note that the special vertices can also be divided into two types - special and hyperspecial. Coloring the vertices using the colors induced by the coloring on $\Theta_{2n}$, the special vertices have colors $0$ (hyperspecial) or $n$ (special but not hyperspecial), while the non-special vertices have two colors each (one color induced from each original vertex): $i$ and $2n-i$ for $1\leq i \less n$. In fact, a vertex is special iff it has color $0$ or $n$ \cite{S09}.

\begin{claim}\label{colors in chambers}
Every chamber in $\Delta_{n}$ has one vertex of color $0$, one vertex of color $n$, and one vertex of color $i,2n-i$ for each $1\leq i \less n$.
\end{claim}

For the proof of this claim see \cite{S09}.

\begin{remark}
We can conclude from \ref{colors in chambers} that $\Delta_{n}$ is an $n$-dimensional simplicial complex.
\end{remark}

\begin{claim}
\label{numOrbits}
The action of $GSp_{2n}(\F)$ on the vertices of $\Delta_{n}$ has $1+\lceil(n-1)/2\rceil$ orbits; One orbit contains all the special vertices and the rest of the orbits contain non-special vertices.
\end{claim}

\begin{proof}
First, the vertex $v=[0,...,0;0,...,0]=[\mathcal{O}^{2n}]$ is a special vertex, as the matrix associated to it is the identity matrix $I_{2n}$ which is obviously in $GSp_{2n}(\F)$. Let $A\in GSp_{2n}(\F)$. From \ref{val and det in GSp}, $\nu (det(A))$ is either equal to $0\mod{2n}$ or $n\mod{2n}$, therefore the color of the vertex $Av$ is $0$ or $n$. Thus, $Av$ is a special vertex. Now let $v'$ be any special vertex in $\Delta_{n}$. Let $A'$ be a matrix associated with $v'$, so that $A'\mathcal{O}^{2n}\in v'$. Thus, $A'v=v'$. And $A'$ is in $GSp_{2n}(\F)$ because $v'$ is a special vertex (i.e. its representatives are symplectic lattices).

Now let $v$ be a non-special vertex of colors $i,2n-i$ for $1\leq i \less n$. Let $A\in GSp_{2n}(\F)$, so $\nu (det(A))$ is either equal to $0\mod{2n}$ or $n\mod{2n}$. In the first case, $Av$ has color $n+i,n-i$, so again it is a non-special vertex, because for $j=n-i$, $1 \leq j \less n$ and $2n-j=n+i$. In the second case, $A$ preserves the color of $v$. Thus, in the orbit of $v$ we can only get vertices of the colors $i,2n-i$ or $n+i,n-i$ for $1\leq i \less n$. Now we can use the fact that $GSp_{2n}(\F)$ is transitive on the chambers\cite{S05}, to conclude that the orbit of $v$ under the $GSp_{2n}(\F)$-action is the set of all the non-special vertices of colors $i,2n-i$ and $n+i,n-i$.
\end{proof}

The next claim has a proof in \cite{S05}:
\begin{claim}
\label{specialVert}
$v=[a_{1},...,a_{n};b_{1},...b_{n}]$ is a special vertex, iff there exists a constant $c\in\Z$ such that for every $1\leq i \leq n$, $a_{i}+b_{i}=c$.
\end{claim}

%% file: Chapter3/ApartmentWeyl.tex
Given a basis for $V$ as in \ref{sec:BTB}, we can define a frame, which is a set of $n$ pairs of lines in $V$: $\{\F u_{i},\F w_{i}\}$ for $1\leq i \leq n$. 
A frame defines an apartment $A$ in $\Delta_{n}$, such that every vertex in $A$ can be written as $[a_{1},...,a_{n};b_{1},...,b_{n}]$ \cite{buildings2}. 
Furthermore, for every apartment $A$ there is a basis for $V$ of hyperbolic pairs such that every vertex in $A$ can be written as $[a_{1},...,a_{n};b_{1},...,b_{n}]$ in that basis. For now, let's fix the basis in \ref{sec:BTB}. 
To build the apartment concretely, we can start by building a chamber containing the vertex $[0,...,0;0,....,0]$. The following vertices form a chamber $C$:
$$[0,...,0;0,....,0],[1,...,1;0,...0],[1,...,1;1,0,...,0]+[0,1,...,1;0,...,0],...,[1,...,1;1,...,1,0]+[0,...0,1;0,...,0].$$ \label{Cn}

Now, to build the rest of the apartment, we can start applying the generators of the affine Weyl group on the chamber $C$. We shall use the following generators (basis vectors not explicitely noted stay fixed):
\label{genWn}
$s_{1}$ swaps between $u_{n}$ and $w_{n}$. Thus, $s_{1}([a_{1},...,a_{n};b_{1},...b_{n}])=[a_{1},...,a_{n-1},b_{n};b_{1},...,b_{n-1},a_{n}]$.

$s_{j}$ for $2\leq j \leq n$ swaps between $u_{n-j+1}$ and $u_{n-j+2}$, and also swaps between $w_{n-j+1}$ and $w_{n-j+2}$. 
Thus, $$s_{j}([a_{1},...,a_{n};b_{1},...b_{n}])=[a_{1},...,a_{n-j},a_{n-j+2},a_{n-j+1},...,a_{n};b_{1},...,b_{n-j},b_{n-j+2},b_{n-j+1},...,b_{n}].$$

$s_{n+1}$ takes $u_{1}$ to $\varpi w_{1}$ and takes $w_{1}$ to $\varpi^{-1} u_{1}$. Thus, $s_{n+1}([a_{1},...,a_{n};b_{1},...b_{n}])=[b_{1}-1,a_{2},...,a_{n};a_{1}+1,b_{2},...,b{n}]$.

These generators satisfy the relations in the Coxeter diagram of $\Tilde{C}_{n}$:
\dynkin[extended,Coxeter]{C}{}
Finally, we have $t$ \cite{M88} which takes $u_{i}$ to $w_{i}$ for $1\leq i \leq n-1$, takes $u_{n}$ to $\varpi w_{n}$, takes $w_{i}$ for $1\leq i \leq n-1$ to $\varpi u_{i}$, and takes $w_{n}$ to $u_{n}$.

\begin{example*}
\label{genW2}
For $n=2$, the generators have the following matrix forms (signs appear so that the matrices belong to $GSp_{4}(\F)$ but of course they don't affect the vertices):
$$s_{1}=\begin{pmatrix}
-1 & 0 & 0 & 0 \\
0  & 0 & 0 & 1 \\
0  & 0 & 1 & 0 \\
0  & 1 & 0 &0 
\end{pmatrix}, \quad s_{2}=\begin{pmatrix}
0 & 1 & 0 & 0 \\
1 & 0 & 0 & 0 \\
0 & 0 & 0 & 1 \\
0 & 0 & 1 & 0 
\end{pmatrix}, \quad
s_{3}=\begin{pmatrix}
0       & 0 & \varpi^{-1} & 0 \\
0       & 1 & 0           & 0 \\
-\varpi & 0 & 0           & 0 \\
0       & 0 & 0           & 1 
\end{pmatrix}, \quad
t=\begin{pmatrix}
0 & 0      & \varpi & 0 \\
0 & 0      & 0      & 1 \\
1 & 0      & 0      & 0 \\
0 & \varpi & 0      & 0 
\end{pmatrix}.$$
By applying words of these generators to the chamber $C$, we get the apartment $A$:
\input{Chapter3/ApartmentDrawing}
\label{apt2}
In the diagram, vertices are labeled by their homothety classes, and edges are solid lines. 
Reflections and the generator $t$ (a rotation) are denoted by blue arrows and labeled with the corresponding (reduced) word in the generators. 
The chamber $C$ is also labeled in blue. The dotted lines on the edges of the diagram indicate that this apartment spreads infinitely in all directions, as it is a tessellation of $\R^2$.
\end{example*}

%% file: Chapter3/ApartmentDrawing.tex
\begin{tikzpicture}
  $$
  \matrix (m) [matrix of math nodes,row sep=2.5cm,column sep=0.04cm]
  {
    \verb|[2,0;0,2]| & \verb|[2,1;0,2]+[2,0;0,1]| & \verb|[2,1;0,1]| & \verb|[2,2;0,1]+[2,1;0,0]| & \verb|[2,2;0,0]| \\
    \verb|[2,0;1,2]+[1,0;0,2]| & \verb|[1,0;0,1]| & \verb|[1,1;0,1]+[1,0;0,0]| & \verb|[1,1;0,0]| & \verb|[2,2;1,0]+[1,2;0,0]| \\
    \verb|[1,0;1,2]| & \verb|[1,0;1,1]+[0,0;0,1]| & \verb|[0,0;0,0]| & \verb|[1,1;1,0]+[0,1;0,0]| & \verb|[1,2;1,0]| \\
    \verb|[1,0;2,2]+[0,0;1,2]| & \verb|[0,0;1,1]| & \verb|[0,1;1,1]+[0,0;1,0]| & \verb|[0,1;1,0]| & \verb|[1,2;2,0]+[0,2;1,0]| \\
    \verb|[0,0;2,2]| & \verb|[0,1;2,2]+[0,0;2,1]| & \verb|[0,1;2,1]| & \verb|[0,2;2,1]+[0,1;2,0]| & \verb|[0,2;2,0]| \\};
  \path
    (m-1-1) edge (m-2-1)
    (m-2-1) edge (m-3-1)
    (m-3-1) edge (m-4-1)
    (m-4-1) edge (m-5-1)
    (m-1-2) edge (m-2-2)
    (m-2-2) edge (m-3-2)
    (m-3-2) edge (m-4-2)
    (m-4-2) edge (m-5-2)
    (m-1-3) edge (m-2-3)
    (m-2-3) edge (m-3-3)
    (m-3-3) edge (m-4-3)
    (m-4-3) edge (m-5-3)
    (m-1-4) edge (m-2-4)
    (m-2-4) edge (m-3-4)
    (m-3-4) edge (m-4-4)
    (m-4-4) edge (m-5-4)
    (m-1-5) edge (m-2-5)
    (m-2-5) edge (m-3-5)
    (m-3-5) edge (m-4-5)
    (m-4-5) edge (m-5-5)
    (m-1-1) edge (m-1-2)
    (m-1-2) edge (m-1-3)
    (m-1-3) edge (m-1-4)
    (m-1-4) edge (m-1-5)
    (m-2-1) edge (m-2-2)
    (m-2-2) edge (m-2-3)
    (m-2-3) edge (m-2-4)
    (m-2-4) edge (m-2-5) 
    (m-3-1) edge (m-3-2)
    (m-3-2) edge (m-3-3)
    (m-3-3) edge (m-3-4)
    (m-3-4) edge (m-3-5)  
    (m-4-1) edge (m-4-2)
    (m-4-2) edge (m-4-3)
    (m-4-3) edge (m-4-4)
    (m-4-4) edge (m-4-5)   
    (m-5-1) edge (m-5-2)
    (m-5-2) edge (m-5-3)
    (m-5-3) edge (m-5-4)
    (m-5-4) edge (m-5-5)   
    (m-1-1) edge (m-2-2)
    (m-2-2) edge (m-3-3)
    (m-3-3) edge (m-4-4)
    (m-4-4) edge (m-5-5)
    (m-1-3) edge (m-2-4)
    (m-2-4) edge (m-3-5)
    (m-3-1) edge (m-4-2)
    (m-4-2) edge (m-5-3)
    (m-5-1) edge (m-4-2)
    (m-4-2) edge (m-3-3)
    (m-3-3) edge (m-2-4)
    (m-2-4) edge (m-1-5)
    (m-1-3) edge (m-2-2)
    (m-2-2) edge (m-3-1)
    (m-3-5) edge (m-4-4)
    (m-4-4) edge (m-5-3);
    \draw [<->,blue] (-0.3,-1.7) -- (0.3,-1.7)
    node [above,right] {s1};
    \draw [<->,blue] (-1.5,-1.65) -- (-2.1,-1.2)
    node [above,right] {s2};
    \draw [<->,blue] (-2.3,-2.7) -- (-2.3,-3.35)
    node [below,right] {s3};
    \draw [<->,blue] (1.5,-1.65) -- (2.1,-1.2)
    node [above,right] {s1s2s1};
    \draw [<->,blue] (1.5,-0.2) -- (1.5,0.2)
    node [above] {s2s1s2};
    \draw [<->,blue] (-2.3,0.8) to [bend left=90] node [below] {t} (-1.3,1.8);
    \draw (2.5,0.7) node[blue,font=\fontsize{20}{10}\sffamily\bfseries]{C};
    \draw [dotted] (8.7,0) -- (9.5,0);
    \draw [dotted] (-8.7,0) -- (-9.5,0);
    \draw [dotted] (0,4.5) -- (0,5.5);
    \draw [dotted] (0,-4.5) -- (0,-5.5);
    \draw [dotted] (8.5,-4.5) -- (9.5,-5);
    \draw [dotted] (8.5,4.5) -- (9.5,5);
    \draw [dotted] (8.5,-4.5) -- (9.5,-5);
    \draw [dotted] (-8.5,-4.5) -- (-9.5,-5);
    \draw [dotted] (-8.5,4.5) -- (-9.5,5);
    $$
\end{tikzpicture}

%% file: Chapter3/IntroCartan.tex
Let $G=GL_n(\F)$ and $K=GL_n(\mathcal{O})$. Denote $$A=\{\begin{psmallmatrix}
\varpi^{m_1} \\
& \varpi^{m_2} \\
& & \cdot \\
& & & \cdot \\
& & & & \cdot \\
& & & & & \varpi^{m_n} \\
\end{psmallmatrix} | m_i \in \Z, m_1 \geq m_2 \geq ... \geq m_n \}$$

\begin{claim} [Cartan Decomposition]
$G=\bigsqcup_{a \in A} KaK$.
\end{claim}

By the theory of symplectic divisors, there is a corresponding decomposition for $G'=GSp_{2n}(\F)$:
$$G'=\bigsqcup_{\omega \in A'}\Gamma\omega\Gamma$$ where $\Gamma=GSp_{2n}(\mathcal{O})$ and $A'$ consists of matrices of the shape $\begin{psmallmatrix}
\varpi^{a_1} \\
&  \cdot \\
& & \cdot \\
& & & \cdot \\
& & & & \varpi^{a_n} \\
& & & & & \varpi^{b_1} \\
& & & & & & \cdot \\
& & & & & & & \cdot \\
& & & & & & & & \cdot \\
& & & & & & & & & \varpi^{b_n} \\
\end{psmallmatrix}$ such that $a_i,b_i\in \Z$, $a_1\leq ... \leq a_n$ and $a_i+b_i=const$.

\begin{remark}
Let $\xi$ be the vertex $[0,...,0;0,...,0]$. Its stabilizer in $G'$ is $\Gamma$ (see proof below). Denote the set of vertices of $\Delta_n$ by $\Delta_n^0$. By the Cartan decomposition for $G'$, $\Delta_n^0=\bigsqcup_{a\in A'}\Gamma a \xi$.
\end{remark}

%% file: Chapter4/IntroGeoOp2.tex
In this chapter we will introduce some operators on the simplices of $\Delta_{2}$. First we will introduce operators on the chambers of $\Delta_{2}$, then we will show how to define operators from vertices to chambers and vice versa, and finally we will combine all these operators to define operators on the vertices of $\Delta_{2}$. We will show that these operators are geometric branching operators, and that the SRW on the vertices of $\Delta_{2}$ can be described by these operators, similarly to the work done in \cite{CP22}. We will start this section by some notations, then proceed in the next section to describe the operators and prove their properties. In the final section we will describe the SRW using these operators.

Denote the vertices of $\Delta_{2}$ by $\Delta_{2}^{0}$. Let $\xi \in \Delta_{2}^{0}$ be the vertex represented by $[0,0;0,0]$ (see definition for this notation in \ref{homothetyExponents}). Note that by \ref{specialVert} $\xi$ is a special vertex.

\begin{remark}
For a general $n$, $\xi=[0,...,0;0,...,0]$. The stabilizer of $\xi$ under the action of $GSp_{2n}(\F)$, which we will denote as $P_{\xi}$, is $GSp_{2n}(\mathcal{O})$; Clearly, $GSp_{2n}(\mathcal{O})\subseteq P_{\xi}$ because for any matrix $M\in GSp_{2n}(\mathcal{O})$, the columns of $M$ span $\mathcal{O}^{2n}$, thus $M\xi=\xi$. In the other direction, for $M\in P_{\xi}$, $M$'s columns span $\mathcal{O}^{2n}$, thus $M$ is in $GSp_{2n}(\mathcal{O})$.
\end{remark}

Let $C$ be the chamber as in example \ref{C2}. Let $\chi$ be the single non-special vertex in $C$, and denote its stabilizer by $P_{\chi}$. Let $S$ be the sector of the apartment $A$ in $\Delta_{2}$ ($A$ as in \ref{apt2}), such that $S$ contains $C$ and $S$ can be identified with the quotient of $\Delta_{2}$ by $P_{\xi}$. More concretely, we choose $S$ to be the following sector: 

\input{Chapter4/drawingS}

\begin{remark}
$S$ is also a fundamental domain for the action of the spherical Weyl group (of type $C_2$) on the building. 
\end{remark}

We identify the special vertices in $S$ with the set $X=\{(x,y)\in \N^{2}|x \geq y\}$ as following: the special vertices in $S$ are of the shape $[x,x+y;y,0]$ such that $x,y\in \N$ and $x \geq y$. This characterization invites the following identification: we identify such a vertex with the pair $(x,y)\in X$. In the final section we will show that the SRW on $X$, projected to $S$ by this identification, can be fully described using the operators in the following section, because every vertex in $S$ can be reached by the operators.

%% file: Chapter4/drawingS.tex
\begin{tikzpicture}
  $$
  \matrix (m) [matrix of math nodes,row sep=1.5cm,column sep=0cm]
  {
    &[0.9cm] &[0.05cm] &[0.05cm] &[0.05cm] \verb|[4,4;0,0]| \\
    & & & \verb|[3,3;0,0]| & \verb|[3,4;0,0]+[4,4;1,0]|\\
    & & \verb|[2,2;0,0]| & \verb|[2,3;0,0]+[3,3;1,0]| & \verb|[3,4;1,0]| \\
    & \verb|[1,1;0,0]| & \verb|[2,2;1,0]+[1,2;0,0]| & \verb|[2,3;1,0]| & \verb|[2,4;1,0]+[3,4;2,0]| \\
    \xi=\verb|[0,0;0,0]| &  \chi=\verb|[1,1;1,0]+[0,1;0,0]| & \verb|[1,2;1,0]| & \verb|[2,3;2,0]+[1,3;1,0]| & \verb|[2,4;2,0]| \\
    };
  \path
    (m-2-4) edge (m-2-5)
    (m-3-3) edge (m-3-4)
    (m-3-4) edge (m-3-5)
    (m-4-2) edge (m-4-3)
    (m-4-3) edge (m-4-4)
    (m-4-4) edge (m-4-5)
    (m-5-1) edge (m-5-2)
    (m-5-2) edge (m-5-3)
    (m-5-3) edge (m-5-4)
    (m-5-4) edge (m-5-5)
    (m-4-2) edge (m-5-2)
    (m-3-3) edge (m-4-3)
    (m-4-3) edge (m-5-3)
    (m-2-4) edge (m-3-4)
    (m-3-4) edge (m-4-4)
    (m-4-4) edge (m-5-4)
    (m-1-5) edge (m-2-5)
    (m-2-5) edge (m-3-5)
    (m-3-5) edge (m-4-5)
    (m-4-5) edge (m-5-5)
    (m-5-1) edge (m-4-2)
    (m-4-2) edge (m-3-3)
    (m-3-3) edge (m-2-4)
    (m-2-4) edge (m-1-5)
    (m-5-3) edge (m-4-4)
    (m-4-4) edge (m-3-5)
    (m-4-2) edge (m-5-3)
    (m-3-3) edge (m-4-4)
    (m-4-4) edge (m-5-5)
    (m-2-4) edge (m-3-5);
    \draw (-5.4,-3.4) node[font=\fontsize{20}{10}\sffamily\bfseries]{C};
    \draw [dotted] (8.7,-4.05) -- (9.5,-4.05);
    \draw [dotted] (8.2,4.5) -- (9,5);
    $$
\end{tikzpicture}

%% file: Chapter4/DescribeOperators2.tex
We will begin by defining geometric branching operators on the chambers of $\Delta_{2}$, which we will denote by $\Delta_{2}^{2}$. $G=GSp_{4}(\F)$ is transitive on $\Delta_{2}^{2}$, thus we can define operators as in the first part of \ref{IntroOperators}, which treats the case of a transitive group action. We will take $C$ to be our $x_{0}$, and denote its stabilizer by $I$.

Note that $I$ stabilizes the chamber $C$ pointwise. This means that it stabilizes both special vertices, but has elements which swap the two homothety classes representing the non-special vertices. By \cite{SCH04}, $I=\{ \begin{psmallmatrix}
\Op & \varpi\Op & \Op & \Op \\
\Op & \Op & \Op & \Op \\
\varpi\Op & \varpi\Op & \Op & \Op \\
\varpi\Op & \varpi\Op & \varpi\Op & \Op\end{psmallmatrix}
\in GSp_4(\Op) \}$. (This formula can be achieved by noticing that changing the basis of the symplectic form $J$ chosen in \cite{SCH04} with the matrix 
$\begin{psmallmatrix}
 & 1 &  &  \\
1 &  &  &  \\
 &  & 1 &  \\
 &  &  & 1 \end{psmallmatrix}$ results in $\Omega^{t}$, thus changing the basis of the transpose of the Iwahori subgroup $I$ in the same paper results in $I$ as written above. It can also be calculated explicitly by calculating intersection of all the stabilizers of the vertices of $C$.)

Let $T_{diag}:\Delta_{2}^{2}\longrightarrow2^{\Delta_{2}^{2}}$ be the operator defined by:
$T_{diag}(C)=Is_{1}ItIs_{2}Is_{1}Is_{2}IC=Is_{1}ts_{2}s_{1}s_{2}C$ for $s_{1},s_{2}$ and $t$ defined in \ref{genW2}, where the final equality is given by the Bruhat relations \cite{buildings2}, and by the fact that $I$ stabilizes $C$. 

Explicitly, 
$T_{diag}(C)=I\begin{psmallmatrix}
\varpi & 0      & 0 & 0 \\
0      & \varpi & 0 & 0 \\
0      & 0      & 1 & 0 \\
0      & 0      & 0 & 1
\end{psmallmatrix}C$. And for every other chamber $C'$ there exists $g\in G$ such that $C'=gC$, thus $$T_{diag}(C')=T_{diag}(gC)=gI\begin{psmallmatrix}
\varpi & 0      & 0 & 0 \\
0      & \varpi & 0 & 0 \\
0      & 0      & 1 & 0 \\
0      & 0      & 0 & 1
\end{psmallmatrix}C.$$ By \ref{IntroOperators} $T_{diag}$ is a geometric branching operator on the chambers of $\Delta_{2}$. 

\begin{claim}
$T_{diag}$ is $\varpi^{3}$-out-regular.
\end{claim}
\label{TdiagReg}
\begin{proof}
Let $s$ be a reflection through one of the walls in $C$, which we will denote by $A$. 
$IsIC$ is a set containing all the chambers which share the wall $A$ with $C$ (not including $C$). By \cite{S09}, there are $\varpi$ such chambers. 
By the Bruhat relations, as $l(Is_{1}ItIs_{2}Is_{1}Is_{2})=l(Is_{1}ts_{2}s_{1}s_{2})=l(Is_{1}) + l(t) + l(s_{2}s_{1}s_{2})= 1+1+1=3$ (we remind that the length of a conjugated reflection is $1$),  $|T_{diag}(C)|=\varpi^{3}$ by \cite{buildings2} (section 6.2). Thus $T_{diag}$ is $\varpi^{3}$-out-regular.
\end{proof}

\begin{remark}
\label{latticesTdiag}
There is another way to show that $T_{diag}$ is $\varpi^{3}$-out-regular, by counting all the vertices which $\xi$ is sent to under $T_{diag}$. We will show it briefly here. 
First, $T_{diag}$ sends $\xi$ to all the vertices whose associated matrices have the Cartan decomposition:
$\begin{psmallmatrix}
\varpi & 0      & 0 & 0 \\
0      & \varpi & 0 & 0 \\
0      & 0      & 1 & 0 \\
0      & 0      & 0 & 1
\end{psmallmatrix}.$
These matrices have the form 
$\begin{psmallmatrix}
\varpi & 0      & a & b \\
0      & \varpi & b & c \\
0      & 0      & 1 & 0 \\
0      & 0      & 0 & 1
\end{psmallmatrix}$, where $a,b,c\in \mathcal{O}/\varpi \mathcal{O}$. There are clearly $\varpi^{3}$ such different matrices, thus $T_{diag}$ is $\varpi^{3}$-out-regular. Furthermore, $GSp_4(\Op)\begin{psmallmatrix}
\varpi & 0      & 0 & 0 \\
0      & \varpi & 0 & 0 \\
0      & 0      & 1 & 0 \\
0      & 0      & 0 & 1
\end{psmallmatrix}I=GSp_4(\Op)\begin{psmallmatrix}
\varpi & 0      & a & b \\
0      & \varpi & b & c \\
0      & 0      & 1 & 0 \\
0      & 0      & 0 & 1
\end{psmallmatrix}=GSp_4(\Op)\begin{psmallmatrix}
\varpi & 0      & 0 & 0 \\
0      & \varpi & 0 & 0 \\
0      & 0      & 1 & 0 \\
0      & 0      & 0 & 1
\end{psmallmatrix}$ by the explicit form of $I$ above and the fact that $I\leq GSp_4(\Op)$. The last equality will be explained in more detail in \ref{T'}.
\end{remark}

\begin{claim}
$T_{diag}$ is $\varpi^{3}$-regular.
\end{claim}
\begin{proof}
This follows from the fact that $G$ acts transitively on $\Delta_{2}^{2}$. Another way to see this is to note that the inverse of $T_{diag}$ is given by $T_{diag}^{-1}(C)=I\begin{psmallmatrix}
1 & 0 & 0      & 0 \\
0 & 1 & 0      & 0 \\
0 & 0 & \varpi & 0 \\
0 & 0 & 0      & \varpi
\end{psmallmatrix}C$ thus a similar analysis shows that it is also $\varpi^{3}$-out-regular.
\end{proof}

\begin{claim}
$T_{diag}$ is collision-free.
\end{claim}
\begin{proof}
There are again (at least) two ways to show this - the first is to note that the vertices to which $T_{diag}$ send $\xi$ are non-equivalent lattices, thus $T_{diag}$ sends $C$ to $\varpi^{3}$ different chambers, which means that it is collision-free. The second way is by showing that $\omega = s_{1}ts_{1}s_{2}s_{1}$ is a translation element, therefore by theorem 5.6 in \cite{LLP20}, $T_{diag}=T_{\omega}$ is collision free. 
To end the proof, $\omega$ is a translation element in the extended affine Weyl group, as it is a translation by the co-root $\alpha_{1}+\alpha_{2}$, with $\alpha_{1},\alpha_{2}$ as defined in \ref{roots coroots C2}.
\end{proof}

Next, we will define one more operator on $\Delta_{2}^{2}$. Let $T_{down}:\Delta_{2}^{2}\longrightarrow2^{\Delta_{2}^{2}}$ be the operator defined by: $T_{down}(C)=Is_{2}Is_{1}Is_{2}Is_{1}ItIs_{1}C=Is_{2}s_{1}s_{2}s_{1}ts_{1}C$, where the final equality is given by the Bruhat relations. Explicitly, $$T_{down}(C)=I\begin{psmallmatrix}
-1 & 0 & 0      & 0      \\
0  & 0 & 0      & \varpi \\
0  & 0 & \varpi & 0      \\
0  & 1 & 0      & 0      \\
\end{psmallmatrix}C.$$ Again, this defines a geometric branching operator on $\Delta_{2}^{2}$.

\begin{claim}
$T_{down}$ is $\varpi^{2}$-regular.
\end{claim}
\begin{proof}
For $\omega=s_{2}s_{1}s_{2}s_{1}ts_{1}$, $l(\omega)=2$, thus $T_{down}$ is $\varpi^{2}$-out-regular. And once again using the fact that $G$ is transitive on $\Delta_{2}^{2}$, we can deduce that $T_{down}$ is $\varpi^{2}$-regular.
\end{proof}

\begin{remark}
\label{latticesTdown}
Once again we can analyze $T_{down}$ by checking which vertices it sends $\xi$ to. The matrices associated to these vertices are of the shape: 
$$\begin{psmallmatrix}
-1 & 0 & 0      & 0      \\
0  & a & 0      & \varpi \\
b  & 0 & \varpi & 0      \\
0  & 1 & 0      & 0      \\
\end{psmallmatrix},\quad a,b\in \mathcal{O}/\varpi \mathcal{O}.$$ We can again see that they represent $\varpi^{2}$ non-equivalent lattices, therefore $T_{down}$ sends $C$ to $\varpi^{2}$ different chambers.
\end{remark}

\begin{claim}
$T_{down}$ is collision-free.
\end{claim}
\begin{proof}
By the above remark we can already deduce that $T_{down}$ is collision-free. 
We will show this again in one more way: we can note that $T_{down}^{2}=I\omega^{2}=
I\begin{psmallmatrix}
1 & 0      & 0          & 0 \\
0 & \varpi & 0          & 0 \\
0 & 0      & \varpi^{2} & 0 \\
0 & 0      & 0          & \varpi
\end{psmallmatrix}$ and $\omega^{2}$ is the translation element in the extended affine Weyl group given by $-\alpha_{1}-2\alpha_{2}$ as defined in \ref{roots coroots C2}. 
Also, $l(\omega^{2})=4=2l(\omega)$, thus again by theorem 5.6 in \cite{LLP20}, $T_{down}$ is collision-free.
\end{proof}

\begin{example*}
The following diagram illustrates the action of $T_{diag}$ and $T_{down}$ on the chamber $C$, restricted to the apartment $A$.

\input{Chapter4/drawingTs}
\end{example*}

To end this section, we will define operators which will take us from vertices of $\Delta_{2}$ to chambers and vice versa. We will use these operators to build geometric branching regular and collision-free operators on $\Delta_{2}^{0}$ in the next section.

According to 3.7 in \cite{CP22}, if $\sigma_{1},\sigma_{2}$ are two simplices (not necessarily of the same dimension) in the building with corresponding stabilizers $P_{1},P_{2}$ respectively, any double coset $P_{1}gP_{2}$ defines a geometric operator from the orbit of $\sigma_{1}$ to subsets of the orbit of $\sigma_{2}$ by $T_{12}(g'\sigma_{1})=g'P_{1}g\sigma_{2}.$ Take $\xi, \chi, C$ as before and denote their stabilizers by $P_{\xi},P_{\chi},I$ respectively.
By the above we can define the following geometric branching operators:
$$T_{0,2}:G\xi\longrightarrow 2^{\Delta_{2}^{2}} \quad \quad T_{0,2}(\xi)=P_{\xi}C=GSp_{4}(\mathcal{O})C$$
$$T_{2,0}:\Delta_{2}^{2}\longrightarrow\Delta_{2}^{0} \quad \quad 
T_{2,0}(C)=I\xi=\xi$$
$$T_{2,0}':\Delta_{2}^{2}\longrightarrow\Delta_{2}^{0} \quad \quad 
T_{2,0}'(C)=I\chi=\chi$$
where the last equalities in the definitions are because $I\leq P_{\xi}$ and $I\leq P_{\chi}$. It is clear that $T_{2,0},T_{2,0}'$ are $1$-out-regular. $T_{0,2}$ is $\prod_{m=1}^{2}\frac{\varpi^{2m}-1}{\varpi-1}$-out-regular by proposition 1.3.5 of \cite{S06}.

%% file: Chapter4/drawingTs.tex
\begin{tikzpicture}
  $$
  \matrix (m) [matrix of math nodes,row sep=2.5cm,column sep=0.04cm]
  {
    \verb|[2,0;0,2]| & \verb|[2,1;0,2]+[2,0;0,1]| & \verb|[2,1;0,1]| & \verb|[2,2;0,1]+[2,1;0,0]| & \verb|[2,2;0,0]| \\
    \verb|[2,0;1,2]+[1,0;0,2]| & \verb|[1,0;0,1]| & \verb|[1,1;0,1]+[1,0;0,0]| & \verb|[1,1;0,0]| & \verb|[2,2;1,0]+[1,2;0,0]| \\
    \verb|[1,0;1,2]| & \verb|[1,0;1,1]+[0,0;0,1]| & \verb|[0,0;0,0]| & \verb|[1,1;1,0]+[0,1;0,0]| & \verb|[1,2;1,0]| \\
    \verb|[1,0;2,2]+[0,0;1,2]| & \verb|[0,0;1,1]| & \verb|[0,1;1,1]+[0,0;1,0]| & \verb|[0,1;1,0]| & \verb|[1,2;2,0]+[0,2;1,0]| \\
    };
  \path
    (m-1-1) edge (m-2-1)
    (m-2-1) edge (m-3-1)
    (m-3-1) edge (m-4-1)
    (m-1-2) edge (m-2-2)
    (m-2-2) edge (m-3-2)
    (m-3-2) edge (m-4-2)
    (m-1-3) edge (m-2-3)
    (m-2-3) edge (m-3-3)
    (m-3-3) edge (m-4-3)
    (m-1-4) edge (m-2-4)
    (m-2-4) edge (m-3-4)
    (m-3-4) edge (m-4-4)
    (m-1-5) edge (m-2-5)
    (m-2-5) edge (m-3-5)
    (m-3-5) edge (m-4-5)
    (m-1-1) edge (m-1-2)
    (m-1-2) edge (m-1-3)
    (m-1-3) edge (m-1-4)
    (m-1-4) edge (m-1-5)
    (m-2-1) edge (m-2-2)
    (m-2-2) edge (m-2-3)
    (m-2-3) edge (m-2-4)
    (m-2-4) edge (m-2-5) 
    (m-3-1) edge (m-3-2)
    (m-3-2) edge (m-3-3)
    (m-3-3) edge (m-3-4)
    (m-3-4) edge (m-3-5)  
    (m-4-1) edge (m-4-2)
    (m-4-2) edge (m-4-3)
    (m-4-3) edge (m-4-4)
    (m-4-4) edge (m-4-5)   
    (m-1-1) edge (m-2-2)
    (m-2-2) edge (m-3-3)
    (m-3-3) edge (m-4-4)
    (m-1-3) edge (m-2-4)
    (m-2-4) edge (m-3-5)
    (m-3-1) edge (m-4-2)
    (m-4-2) edge (m-3-3)
    (m-3-3) edge (m-2-4)
    (m-2-4) edge (m-1-5)
    (m-1-3) edge (m-2-2)
    (m-2-2) edge (m-3-1)
    (m-3-5) edge (m-4-4);
    \draw (6.2,2.2) node[blue,font=\fontsize{15}{10}\sffamily\bfseries]{$T_{diag}$(C)};
    \draw (1.3,-3.7) node[brown,font=\fontsize{15}{10}\sffamily\bfseries]{$T_{down}$(C)};
    \draw (2.5,-0.7) node[black,font=\fontsize{20}{10}\sffamily\bfseries]{C};
    $$
\end{tikzpicture}

%% file: Chapter4/SRW2.tex
We are now going to describe the allowed moves in a SRW on the vertices of $\Delta_{2}$, starting at the vertex $\xi$, using the operators defined in the previous section. 
Denote the projection from $\Delta_{2}^{0}$ to vertices in $A$, given by the Cartan decomposition (i.e. it sends a vertex to the vertex corresponding to the diagonal matrix in its associated matrix's Cartan decomposition), by $\Phi$. 

\begin{claim}
\label{T'}
For any $x\geq y\in\Z_{\geq 0}$, we have $\Phi^{-1}(\begin{psmallmatrix}
\varpi^{x} & 0   & 0          & 0 \\
0 & \varpi^{x+y} & 0          & 0 \\
0 & 0            & \varpi^{y} & 0 \\
0 & 0            & 0          & 1
\end{psmallmatrix})= T_{x,y}'(\xi)$, where
$$T_{x,y}'=T_{2,0}\circ T_{down}^{2y} \circ T_{diag}^{x+y} \circ T_{0,2}.$$
\end{claim}
\begin{proof}
$$T_{x,y}'(\xi)=
P_{\xi}I
(I\begin{psmallmatrix}
\varpi & 0 & 0 & 0 \\
0 & \varpi & 0 & 0 \\
0 & 0      & 1 & 0 \\
0 & 0      & 0 & 1
\end{psmallmatrix}I)^{x+y}
(I\begin{psmallmatrix}
1 & 0 & 0 & 0 \\
0 & \varpi & 0 & 0 \\
0 & 0 & \varpi^{2} & 0 \\
0 & 0 & 0 & \varpi
\end{psmallmatrix}I)^{y}\xi.$$ 
Note first that $I\leq P_{\xi}$, as I stabilizes all the vertices in $C$, and in particular $\xi$.
Second, since the matrices $\begin{psmallmatrix}
1 & 0 & 0 & 0 \\
0 & \varpi & 0 & 0 \\
0 & 0 & \varpi^{2} & 0 \\
0 & 0 & 0 & \varpi
\end{psmallmatrix}$ and $\begin{psmallmatrix}
\varpi & 0 & 0 & 0 \\
0 & \varpi & 0 & 0 \\
0 & 0      & 1 & 0 \\
0 & 0      & 0 & 1
\end{psmallmatrix}$ are translation elements in the affine Weyl group $\mathcal{W}$, and for every translation element $w\in \mathcal{W}$, $(IwI)^j=Iw^jI$, we conclude that $$T_{x,y}'(\xi)=
P_{\xi}
\begin{psmallmatrix}
\varpi^{x+y} & 0 & 0 & 0 \\
0 & \varpi^{x+y} & 0 & 0 \\
0 & 0      & 1 & 0 \\
0 & 0      & 0 & 1
\end{psmallmatrix}I
\begin{psmallmatrix}
1 & 0 & 0 & 0 \\
0 & \varpi^y & 0 & 0 \\
0 & 0 & \varpi^{2y} & 0 \\
0 & 0 & 0 & \varpi^y
\end{psmallmatrix}\xi.$$
By the computations in \ref{latticesTdiag}, 
$$P_{\xi}(I\begin{psmallmatrix}
\varpi & 0 & 0 & 0 \\
0 & \varpi & 0 & 0 \\
0 & 0      & 1 & 0 \\
0 & 0      & 0 & 1
\end{psmallmatrix}I)^{x+y}=P_{\xi}\{\begin{psmallmatrix}
\varpi^{x+y} & 0      & a & b \\
0      & \varpi^{x+y} & b & c \\
0      & 0      & 1 & 0 \\
0      & 0      & 0 & 1
\end{psmallmatrix}|a,b,c\in \mathcal{O}/\varpi^{x+y} \mathcal{O}\}$$
By multiplying from the left with matrices of the shape $\begin{psmallmatrix}
1 & 0 & * & 0 \\
0 & 1 & 0 & 0 \\
0 & 0 & 1 & 0 \\
0 & 0 & 0 & 1
\end{psmallmatrix},
\begin{psmallmatrix}
1 & 0 & 0 & 0 \\
0 & 1 & 0 & * \\
0 & 0 & 1 & 0 \\
0 & 0 & 0 & 1
\end{psmallmatrix},
\begin{psmallmatrix}
1 & 0 & 0 & * \\
0 & 1 & * & 0 \\
0 & 0 & 1 & 0 \\
0 & 0 & 0 & 1
\end{psmallmatrix}$ where $*\in \mathcal{O}$, which are all in $P_{\xi}=GSp_{4}(\mathcal{O})$, we can see that $P_{\xi}\begin{psmallmatrix}
\varpi^{x+y} & 0      & a & b \\
0      & \varpi^{x+y} & b & c \\
0      & 0            & 1 & 0 \\
0      & 0            & 0 & 1
\end{psmallmatrix}=P_{\xi}\begin{psmallmatrix}
\varpi^{x+y} & 0      & 0 & 0 \\
0      & \varpi^{x+y} & 0 & 0 \\
0      & 0            & 1 & 0 \\
0      & 0            & 0 & 1
\end{psmallmatrix}$.

Thus, 
$$T_{x,y}'(\xi)=P_{\xi}\begin{psmallmatrix}
\varpi^{x+y} & 0      & a & b \\
0      & \varpi^{x+y} & b & c \\
0      & 0      & 1 & 0 \\
0      & 0      & 0 & 1
\end{psmallmatrix}
\begin{psmallmatrix}
1 & 0      & 0 & 0 \\
0      & \varpi^{y} & 0 & 0 \\
0      & 0      & \varpi^{2y} & 0 \\
0      & 0      & 0 & \varpi^{y}
\end{psmallmatrix}\xi= P_{\xi}\begin{psmallmatrix}
\varpi^{x+y} & 0      & 0 & 0 \\
0      & \varpi^{x+y} & 0 & 0 \\
0      & 0      & 1 & 0 \\
0      & 0      & 0 & 1
\end{psmallmatrix}
\begin{psmallmatrix}
1 & 0      & 0 & 0 \\
0      & \varpi^{y} & 0 & 0 \\
0      & 0      & \varpi^{2y} & 0 \\
0      & 0      & 0 & \varpi^{y}
\end{psmallmatrix}\xi$$.

Therefore, $$T_{x,y}'(\xi)=P_{\xi}\begin{psmallmatrix}
\varpi^{x+y} & 0  & 0 & 0 \\
0 & \varpi^{x+2y} & 0 & 0 \\
0      & 0      & \varpi^{2y} & 0 \\
0      & 0      & 0 & \varpi^{y}
\end{psmallmatrix}\xi.$$

As $\varpi^{y}\xi=\xi$ (because $\xi$ is a homothety class),
$T_{x,y}'(\xi)=P_{\xi}\begin{psmallmatrix}
\varpi^{x} & 0  & 0 & 0 \\
0 & \varpi^{x+y} & 0 & 0 \\
0      & 0      & \varpi^{y} & 0 \\
0      & 0      & 0 & 1
\end{psmallmatrix}\xi = \Phi^{-1}(\begin{psmallmatrix}
\varpi^{x} & 0   & 0          & 0 \\
0 & \varpi^{x+y} & 0          & 0 \\
0 & 0            & \varpi^{y} & 0 \\
0 & 0            & 0          & 1
\end{psmallmatrix})$.
\end{proof}

\begin{remark}
The matrices associated to the special vertices of the sector $S$ are of the shape 
$$\begin{psmallmatrix}
\varpi^{x} & 0   & 0          & 0 \\
0 & \varpi^{x+y} & 0          & 0 \\
0 & 0            & \varpi^{y} & 0 \\
0 & 0            & 0          & 1
\end{psmallmatrix}.$$ Thus, given any special vertex in $\Delta_{2}$, $T_{x,y}'$ takes $\xi$ to a subset of the special vertices containing the given special vertex.
\end{remark}

\begin{claim}
For any $x,y\in\Z_{\geq 0}$, we have $\Phi^{-1}(\begin{psmallmatrix}
\varpi^{x} & 0   & 0          & 0 \\
0 & \varpi^{x+y+1} & 0          & 0 \\
0 & 0            & \varpi^{y} & 0 \\
0 & 0            & 0          & 1
\end{psmallmatrix})=T_{x,y}''(\xi)$, where
$$T_{x,y}''=T_{2,0}'\circ T_{down}^{2y} \circ T_{diag}^{x+y} \circ T_{0,2}.$$
\end{claim}
\begin{proof}
$T_{x,y}''(\xi)=
P_{\xi}I
(I\begin{psmallmatrix}
\varpi & 0 & 0 & 0 \\
0 & \varpi & 0 & 0 \\
0 & 0      & 1 & 0 \\
0 & 0      & 0 & 1
\end{psmallmatrix}I)^{x+y}
(I\begin{psmallmatrix}
1 & 0 & 0 & 0 \\
0 & \varpi & 0 & 0 \\
0 & 0 & \varpi^{2} & 0 \\
0 & 0 & 0 & \varpi
\end{psmallmatrix}I)^{y}\chi.$ By the proof of \ref{T'},
$$T_{x,y}''(\xi)=P_{\xi}\begin{psmallmatrix}
\varpi^{x} & 0   & 0          & 0 \\
0 & \varpi^{x+y} & 0          & 0 \\
0 & 0            & \varpi^{y} & 0 \\
0 & 0            & 0          & 1
\end{psmallmatrix}\chi=P_{\xi}\begin{psmallmatrix}
\varpi^{x} & 0   & 0          & 0 \\
0 & \varpi^{x+y} & 0          & 0 \\
0 & 0            & \varpi^{y} & 0 \\
0 & 0            & 0          & 1
\end{psmallmatrix}
\begin{psmallmatrix}
1 & 0      & 0 & 0 \\
0 & \varpi & 0 & 0 \\
0 & 0      & 1 & 0 \\
0 & 0      & 0 & 1
\end{psmallmatrix}\xi=P_{\xi}\begin{psmallmatrix}
\varpi^{x} & 0   & 0          & 0 \\
0 & \varpi^{x+y+1} & 0          & 0 \\
0 & 0            & \varpi^{y} & 0 \\
0 & 0            & 0          & 1
\end{psmallmatrix}\xi=\Phi^{-1}(\begin{psmallmatrix}
\varpi^{x} & 0   & 0          & 0 \\
0 & \varpi^{x+y+1} & 0          & 0 \\
0 & 0            & \varpi^{y} & 0 \\
0 & 0            & 0          & 1
\end{psmallmatrix}).$$
\end{proof}

\begin{remark}
Every non-special vertex in the sector $S$ has one of the two associated matrices given by a translation of a special vertex by the matrix $\begin{psmallmatrix}
1 & 0      & 0 & 0 \\
0 & \varpi & 0 & 0 \\
0 & 0      & 1 & 0 \\
0 & 0      & 0 & 1
\end{psmallmatrix}.$ Thus, it is of the shape 
$$\begin{psmallmatrix}
\varpi^{x} & 0   & 0          & 0 \\
0 & \varpi^{x+y+1} & 0          & 0 \\
0 & 0            & \varpi^{y} & 0 \\
0 & 0            & 0          & 1
\end{psmallmatrix}.$$ Therefore, given any non-special vertex in $\Delta_{2}$, $T_{x,y}''$ takes $\xi$ to a subset of the non-special vertices containing the given non-special vertex.
\end{remark}

\begin{corollary}
\label{corollary SRW2}
Define $T_{x,y}:G\xi\longrightarrow2^{\Delta_{2}^{0}}$ by 
$$T_{x,y}(\xi)=T_{x,y}'(\xi)\bigsqcup T_{x,y}''(\xi).$$ By the above, any vertex in $\Delta_{2}^{0}$ is contained in $T_{x,y}(\xi)$ for the appropriate $x,y\in \Z_{\geq 0}$, and $T_{x,y}$ is a geometric branching collision-free operator. Furthermore, every vertex is reached exactly once by $T_{x,y}(\xi)$ when iterating over different values of $x$ and $y$. Therefore, it is possible to use $T_{x,y}$ to understand the rate in which the SRW covers balls around a starting vertex (for example $\xi$), and use it to establish cutoff and calculate the cutoff time in a similar way to the techniques used in \cite{CP22}.
\end{corollary}

%% file: Chapter5/IntroGeoOpn.tex
In this chapter we will generalize the work of the previous chapter to any $n\geq 2$. We will introduce geometric branching operators on the chambers of $\Delta_{n}$, and operators taking vertices to chambers and vice versa.
Combining all these operators together, we will get a family of regular collision-free geometric branching operators which will enable us to described the allowed moves in a SRW on $\Delta_{n}^{0}$, the vertices of $\Delta_{n}$. 

Throughout this chapter, let:

$\xi=[0,...,0;0,...,0]$,

$\chi_{1}=[0,...0,1;0,0,...,0]+[1,...,1;1,...,1,0]$,

$\chi_{2}=[0,...,0,1,1;0,0,...,0]+[1,...,1;1,...,1,0,0]$,

...

$\chi_{n-1}=[0,1,...,1;0,...,0]+[1,...,1;1,0,...,0]$, 

$\xi'=[1,...,1;0,...,0]$. 

By \ref{Cn}, $\xi, \chi_{1},...,\chi_{n-1},\xi'$ are vertices of a chamber. We will denote that chamber as $C$. By \ref{specialVert} it is clear that $\xi, \xi'$ are special vertices and the $\chi_{i}$ are non-special vertices for $1\leq i \leq n-1$. 
For a vertex $\sigma$, denote its stabilizer by $P_{\sigma}$, and denote the stabilizer of $C$ by $I$. Let $S$ be the sector of the apartment $A$ which consists of $C$ and its reflections by the Weil group $\mathcal{W}$, such that $S$ contains $C$ and $S$ can be identified with the quotient of $\Delta_{n}$ by $P_{\xi}$.

\begin{remark}
Similarly to the case for $n=2$, $S$ is a fundamental domain for the action of the spherical Weyl group on the building.
\end{remark}

We identify the special vertices in $S$ with the set $X=\{(x_{1},...,x_{n})\in \N^{n}|x_{1} \geq x_{2} \geq ... \geq x_{n} \}$ as following: the special vertices in $S$ are of the shape $[x_{1}, x_{1}+x_{2},...,x_{1}+...+x_{n};x_{2}+...+x_{n},x_{3}+...+x_{n},...,x_{n},0]$ such that $x_{i}\in \N$ and $x_{i} \geq x_{i+1}$. 
This characterization invites the following identification: we identify such a vertex with  $(x_{1},...,x_{n})\in X$. In the final section we will show that the SRW on $X$, projected to $S$ by this identification, can be fully described using the operators in the following section.

%% file: Chapter5/DescribeOperatorsn.tex
We will again begin by describing operators on the chambers of $\Delta_{n}$, which we denote as $\Delta_{n}^{n}$. $G=GSp_{2n}(K)$ is transitive on $\Delta_{n}^{n}$, thus we can once again use the theory in the first part of \ref{IntroOperators} to define the operators.

Let $T_{diag}:\Delta_{n}^{n}\longrightarrow2^{\Delta_{n}^{n}}$ be the operator defined by:
$T_{diag}(C)=Is_{1}ItIs_{2}Is_{1}Is_{2}C=Is_{1}ts_{2}s_{1}s_{2}C$ for $s_{1},s_{2}$ and $t$ defined in \ref{genWn}, where the final equality is given by the Bruhat relations. Explicitly, $T_{diag}(C)=I\begin{psmallmatrix}
\varpi  \\
& ...   \\
&  & \varpi  \\
&  &  & 1  \\
&  &  &  & ... \\
&  &  &  &  & 1
\end{psmallmatrix}C$ (the matrix is a diagonal matrix, with the first $n$ diagonal entries equal to $\varpi$ and the last $n$ are equal to $1$). Just like in the case of $n=2$, $T_{diag}$ is a $\varpi^3$-regular collision-free geometric branching operator.

Now we will define a family of operators which generalize $T_{down}^2$ from the previous chapter. We will generalize $T_{down}^2$ as we used it in the formula for $T_{x,y}$, and we didn't specifically need $T_{down}$ (not squared) for describing the SRW using the operators. 
We'll start by defining $\sigma_{i}=s_{i}s_{i-1}\cdot \cdot \cdot s_{2}s_{1}s_{2}\cdot \cdot \cdot s_{i-1}s_{i}$ for $1\leq i \leq n$. 

\begin{claim}
$\sigma_{i}$ interchanges $u_{n-i+1}$ and $w_{n-i+1}$, and keeps the other basis vectors in place.
\end{claim}
\begin{proof}
By induction on $i$. For $i=1$ it is the definition of $s_{1}$. Assume that $\sigma_{i-1}$ interchanges $u_{n-(i-1)+1}=u_{n-i+2}$ and $w_{n-(i-1)+1}=w_{n-i+2}$. 
By the definition of $s_{i}$ ($i\geq 2$), it interchanges $u_{n-i+1}$ with $u_{n-i+2}$ and $w_{n-i+1}$ with $w_{n-i+2}$. 
Thus $s_{i}\sigma_{i-1}s_{i}$ sends $u_{n-i+1}$ first to $u_{n-i+2}$, then to $w_{n-i+2}$ and finally to $w_{n-i+1}$. 
Similarly, it sends $w_{n-i+1}$ first to $w_{n-i+2}$, then to $u_{n-i+2}$ and finally to $u_{n-i+1}$. It sends $u_{n-i+2}$ and $w_{n-i+2}$ to themselves (as only $s_{i}$ affects them). And it keeps all the other basis vectors in place.
\end{proof}

Define $t'=s_{1}ts_{1}$. It is an easy exercise to check that $t'$ does the following to the basis vectors:
$$u_{i}\mapsto w_{i}, \quad w_{i}\mapsto \varpi u_{i}.$$ Finally, define $\omega_{i}=\sigma_{n-i+1}\sigma_{n-i+2} \cdot \cdot \cdot \sigma_{n-1} \sigma_{n}t'$.

\begin{claim}
Let $1\leq j \leq n-1$. $\omega_{j}^2$ is the following mapping: $u_{i}\mapsto \varpi u_{i}, w_{i}\mapsto \varpi w_{i}$ for $j \less i \leq n$. $u_{i}\mapsto u_{i}$ and $w_{i}\mapsto \varpi^2 w_{i}$ for $1\leq i \leq j$.
\end{claim}

\begin{proof}
As $\sigma_{n-r+1}$ affects only the basis vectors $u_{r},w_{r}$, only $t'$ affects the rest of the basis vectors. Thus for $i$, $j \less i \leq n$, $\omega_{j}^2$ sends $u_{i}$ to $\varpi u_{i}$ and $w_{i}$ to $\varpi w_{i}$. 
Now, for $1 \leq i \leq j$, $\omega_{j}^2$ does the following to $u_{i}$: first $t'$ sends it to $w_{i}$, then $\sigma_{n-i+1}$ sends it back to $u_{i}$, and the same two steps happens once more (the rest of the $\sigma$s don't affect it). Thus $u_{i}$ is sent to itself. 
Last, $w_{i}$ is sent by $t'$ to $\varpi u_{i}$, which is sent by $\sigma_{n-i+1}$ to $\varpi w_{i}$, which is sent by $t'$ to $\varpi^2 u_{i}$, which is finally sent by $\sigma_{n-i+1}$ to $\varpi^2 w_{i}$.
\end{proof}

\begin{remark}
The matrix forms of the $\omega_{j}^2$ are as following:
$$\omega_{1}^2=\begin{psmallmatrix}
1 \\
& \varpi  \\
& ...   \\
&  & \varpi  \\
&  &  & \varpi^{2}  \\
&  &  &  & \varpi \\
&  &  &  &  & ...\\
&  &  &  &  &  &  \varpi \\
\end{psmallmatrix}$$
$$\omega_{2}^2=\begin{psmallmatrix}
1 \\
& 1   \\
&  & \varpi  \\
&  &  & ... \\
&  &  &  & \varpi  \\
&  &  &  &  & \varpi^2 \\
&  &  &  &  &  & \varpi^2\\
&  &  &  &  &  &  & \varpi\\
&  &  &  &  &  &  & & ...\\
&  &  &  &  &  &  & & & \varpi\\
\end{psmallmatrix}$$
$$...$$
$$\omega_{n-1}^2=\begin{psmallmatrix}
1 \\
& ...   \\
&  & 1  \\
&  &  & \varpi \\
&  &  &  & \varpi^2  \\
&  &  &  &  & ... \\
&  &  &  &  &  & \varpi^2\\
&  &  &  &  &  &  & \varpi\\
\end{psmallmatrix}$$
\end{remark}

\begin{definition}
Define $T_{i}(C)=I\omega_{i}^2C$ for $1\leq i \leq n-1$.
\end{definition}

By \ref{IntroOperators}, the operators $T_{i}$ are geometric branching operators. We can think of them as the generalization of $T_{down}^2$ in \ref{sec:Main proof n=2}. Each one of them is regular with a different regularity (this can be seen by using the length function like in \ref{TdiagReg}, and because $G$ is transitive on $\Delta_{n}^{n}$). They are also collision-free as they are translation elements by coroots, like in the case of $T_{down}^2$. Thus we have defined $n-1$ regular collision-free geometric branching operators.

To end this section, we will once again define operators from vertices to chambers and vice versa. Begin by defining the operator 
$$T_{0,n}:G\xi \longrightarrow 2^{\Delta_{n}^{n}} \quad \quad T_{0,n}(\xi)=P_{\xi}C=GSp_{2n}(\mathcal{O})C.$$ 
Next, we will define, for every orbit of $G$ on the vertices, an operator which will take a chamber and give a vertex in that chamber belonging to the given orbit. 
The first orbit we will choose is the special vertices' orbit. Define $$T_{n,0}:\Delta_{n}^{n}\longrightarrow\Delta_{n}^{0} \quad \quad 
T_{n,0}(C)=I\xi=\xi$$
Next, choose an orbit of non-special vertices of colors $i,2n-i$ and $n+i,n-i$ for some $1 \leq i \leq \lceil \frac{n}{2} \rceil$ (see \ref{numOrbits}). Let $\chi_{i}$ be the single vertex of color $i,2n-i$ in $C$. Define 
$$T_{n,0}^{(i)}:\Delta_{n}^{n}\longrightarrow G\chi_{i} \quad \quad 
T_{n,0}^{(i)}=I\chi_{i}=\chi_{i}$$ where the last equality in the definitions are because $I\leq P_{\xi}$ and $I\leq P_{\chi_{i}}$. It is clear that $T_{n,0},T_{n,0}^{(i)}$ are $1$-out-regular. $T_{0,n}$ is $\prod_{m=1}^{n}\frac{\varpi^{2m}-1}{\varpi-1}$-out-regular by proposition 1.3.5 of \cite{S06}.

%% file: Chapter5/SRWn.tex
We are now going to describe the allowed moves in a SRW on $\Delta_{n}^{0}$, starting at the vertex $\xi$, like in \ref{SRW2}. Denote the projection from $\Delta_{n}^{0}$ to vertices in $A$ given by the Cartan decomposition, by $\Phi$.

\begin{claim}
For any $x_1\geq x_2 \geq ... \geq x_n \in\Z_{\geq 0}$, $\Phi^{-1}(\begin{psmallmatrix}
\varpi^{x_{1}} \\
& \varpi^{x_{1}+x_{2}} \\
& & ... \\
& & & \varpi^{x_{1}+...+x_{n}} \\
& & & & \varpi^{x_{2}+...+x_{n}} \\
& & & & & ... \\
& & & & & & \varpi^{x_{n}} \\
& & & & & & & 1 \\
\end{psmallmatrix})=T'_{x_{1},...,x_{n}}(\xi)$, where 
$$T'_{x_{1},...,x_{n}}=T_{n,0} \circ T_{n-1}^{x_{n}} \circ ... \circ T_{2}^{x_{3}} \circ T_{1}^{x_{2}} \circ T_{diag}^{x_{1}+...+x_{n}} \circ T_{0,n}$$
\end{claim}
\begin{proof}
Using the fact that for every translation element $w\in \mathcal{W}$, $(IwI)^j=Iw^jI$,
$$T'_{x_{1},...,x_{n}}(\xi)=P_{\xi}I
\begin{psmallmatrix}
\varpi^{x_1+...+x_n}  \\
& ...   \\
&  & \varpi^{x_1+...+x_n}  \\
&  &  & 1  \\
&  &  &  & ... \\
&  &  &  &  & 1
\end{psmallmatrix}I
\begin{psmallmatrix}
1  \\
& ...   \\
&  & 1  \\
&  &  & \varpi^{x_n}  \\
&  &  &  & \varpi^{2x_n} \\
&  &  &  &  & ... \\
&  &  &  &  &  & \varpi^{2x_n} \\
&  &  &  &  &  &  & \varpi^{x_n}  \\
\end{psmallmatrix}I ... 
\begin{psmallmatrix}
1  \\
& \varpi^{x_2}   \\
&  & ...\\
&  &  & \varpi^{x_2}  \\
&  &  &  & \varpi^{2x_2} \\
&  &  &  &  &  \varpi^{x_2}\\
&  &  &  &  &  & ... \\
&  &  &  &  &  &  & \varpi^{x_2}  \\
\end{psmallmatrix}\xi.$$
Using similar computations to the case for $n=2$, $I\leq P_{\xi}$, and 
$$T'_{x_{1},...,x_{n}}(\xi)=P_{\xi}
\begin{psmallmatrix}
\varpi^{x_1+...+x_{n-1}}  \\
& ...   \\
&  & \varpi^{x_1+...+x_n}  \\
&  &  & \varpi^{x_n}  \\
&  &  &  & ... \\
&  &  &  &  & \varpi^{x_n} \\
&  &  &  &  &  & 1
\end{psmallmatrix}I
 ... 
\begin{psmallmatrix}
1  \\
& \varpi^{x_2}   \\
&  & ...\\
&  &  & \varpi^{x_2}  \\
&  &  &  & \varpi^{2x_2} \\
&  &  &  &  &  \varpi^{x_2}\\
&  &  &  &  &  & ... \\
&  &  &  &  &  &  & \varpi^{x_2}  \\
\end{psmallmatrix}\xi.$$
Continuing this way, at the final step
$$T'_{x_{1},...,x_{n}}(\xi)=P_{\xi}
\begin{psmallmatrix}
\varpi^{x_1+x_{2}}  \\
& \varpi^{x_1+2x_{2}}   \\
&  & ...  \\
&  & &\varpi^{x_1+2x_{2}+x_3+...+x_n} \\
&  &  &  & \varpi^{2x_{2}+x_3+...+x_n} \\
&  &  &  &  & ... \\
&  &  &  &  &  & \varpi^{x_2+x_n} \\
&  &  &  &  &  & &  \varpi^{x_2} \\
\end{psmallmatrix}\xi=$$

$$=P_{\xi}
\begin{psmallmatrix}
\varpi^{x_1}  \\
& \varpi^{x_1+x_{2}}   \\
&  & ...  \\
&  & &\varpi^{x_1+x_{2}+x_3+...+x_n} \\
&  &  &  & \varpi^{x_{2}+x_3+...+x_n} \\
&  &  &  &  & ... \\
&  &  &  &  &  & \varpi^{x_n} \\
&  &  &  &  &  & &  1 \\
\end{psmallmatrix}\xi=\Phi^{-1}(\begin{psmallmatrix}
\varpi^{x_{1}} \\
& \varpi^{x_{1}+x_{2}} \\
& & ... \\
& & & \varpi^{x_{1}+...+x_{n}} \\
& & & & \varpi^{x_{2}+...+x_{n}} \\
& & & & & ... \\
& & & & & & \varpi^{x_{n}} \\
& & & & & & & 1 \\
\end{psmallmatrix}).$$
\end{proof}

Now, to reach orbits of non-special vertices, we can define for $1\leq i \leq \lceil \frac{n}{2} \rceil$ 
$$T^{(i)}_{x_{1},...,x_{n}}=T_{n,0}^{(i)} \circ T_{n-1}^{x_{n}} \circ ... \circ T_{2}^{x_{3}} \circ T_{1}^{x_{2}} \circ T_{diag}^{x_{1}+...+x_{n}} \circ T_{0,n}$$. 

\begin{claim}
For any $x_1\geq x_2 \geq ... \geq x_n \in\Z_{\geq 0}$, $\Phi^{-1}(\begin{psmallmatrix}
\varpi^{x_{1}} \\
& \varpi^{x_{1}+x_{2}} \\
& & ... \\
& & & \varpi^{x_{1}+...+x_{n}} \\
& & & & \varpi^{x_{2}+...+x_{n}} \\
& & & & & ... \\
& & & & & & \varpi^{x_{n}} \\
& & & & & & & 1 \\
\end{psmallmatrix}
\begin{psmallmatrix}
1 \\
& ... \\
& & \varpi \\
& & & ... \\
& & & & 1 \\
& & & & & ... \\
& & & & & & 1 \\
\end{psmallmatrix})=T^{(i)}_{x_{1},...,x_{n}}(\xi)$
\end{claim}
\begin{proof}
Just like the previous proof.
\end{proof}

Finally, define 
$$T_{x_{1},...,x_{n}}(\xi)=\bigsqcup_{i=1}^{\lceil \frac{n}{2} \rceil}T^{(i)}_{x_{1},...,x_{n}}(\xi)\bigsqcup T'_{x_{1},...,x_{n}}(\xi)$$

\begin{corollary}
$T_{x_{1},...,x_{n}}$ reaches all the vertices in $\Delta_{n}^{0}$ for the appropriate $x_1,...,x_n\in \Z_{\geq 0}$. Additionally, $T_{x_{1},...,x_{n}}$ is a geometric branching and collision-free operator, thus the conclusion in \ref{corollary SRW2} is also valid in the general case.
\end{corollary}

%% file: Discussion.tex
This work leaves some open questions. For example, it is known from previous work that the SRW on Ramanujan complexes arising from the group $Sp_{n}$ over a non-archimedean local field exhibits cutoff. However, what is the cutoff time for those SRWs? This work might be a step in establishing the foundations for answering this question, as done in \cite{CP22}. A more general question is, what is the cutoff time for a SRW on Ramanujan complexes arising from a general algebraic group over a non-archimedean local field?